\documentclass[12pt,reqno,twoside]{amsart}
\usepackage{graphicx} 
\usepackage{amsmath}
\usepackage{amssymb}

\usepackage{color}

\newcommand{\ignore}[1]{}

\newcommand{\Q}{{\mathbb {Q}}}

\newcommand{\R}{{\mathbb{R}}}  
\newcommand{\Z}{{\mathbb{Z}}}

\newcommand{\h}{{\hbox{H}_f}}

\newcommand{\N}{{\mathbb{N}}}

\newcommand\vare{\varepsilon}

\newcommand\NN{\mathbb N}
\newcommand\RR{\mathbb R}
\newcommand\ZZ{\mathbb Z}
\newcommand\QQ{\mathbb Q}

\newcommand\cB{\mathcal{B}}

\newcommand\cW{\mathcal{W}}

\newcommand\cU{\mathcal{U}}

\newcommand\cI{\mathcal{I}}

\newcommand\supp{\operatorname{supp}}

\newcommand\vol{\operatorname{vol}}

\newcommand\norm[1]{\left\|#1\right\|}

\newcommand\abs[1]{\left|#1\right|}

\newcommand\inn[1]{\left\langle #1 \right\rangle}
\newcommand\set[1]{\left\{{#1}\right\}}

\newtheorem{thm}{Theorem}[section]
\newtheorem{lem}[thm]{Lemma}

\newtheorem{cor}[thm]{Corollary}

\newtheorem{rem}[thm]{Remark}

\numberwithin{equation}{section}

\begin{document}

\title[Explicit discrepancy estimates]{Discrepancy of rational points \\ in simple algebraic groups}
\author{Alexander Gorodnik and Amos Nevo} 
\address{Institute f\"ur Mathematik, Universit\"at Z\"urich }
\email{alexander.gorodnik@math.uzh.ch}
\address{Department of Mathematics, Technion IIT, Israel}
\email{anevo@tx.technion.ac.il}


\date{\today}



\begin{abstract} 
The present paper analyzes the discrepancy of distribution of rational points on general semisimple algebraic group varieties. The results include mean-square, almost sure, and uniform discrepancy estimates with explicit error bounds, which apply to general families of subsets, and are valid at arbitrarily small scales. We also consider an analogue of W. Schmidt's classical theorem, which establishes effective almost sure asymptotic counting of rational solutions to Diophantine inequalities in Euclidean spaces. We formulate and prove a version of it for rational points on the group variety, together with an effective bound which in some instances can be expected to be best possible.

\end{abstract}

\maketitle


\section{Introduction}\label{sec:intro}
 
Our goal is analyzing the 
discrepancy of distribution for rational points on certain homogeneous algebraic varieties; 
that is, the behaviour of the 
counting function for 
the number of rational solutions of Diophantine inequalities.
The most classical setting for this problem is that of rational points in Euclidean spaces, and let us begin by recalling some of the most basic results.  
For vectors ${x}\in \mathbb{R}^d$,
one considers the inequality 
\begin{equation}\label{ineq-Khinchin}
\norm{{x} - {{ p}}/{q}}_\infty < {\psi(q)},\quad \hbox{with $({ p},q) \in \Z^d\times \N,$}
\end{equation}
where $\|\cdot\|_\infty$ denotes the maximum norm on $\R^d$, and 
$\psi :(0,\infty)\to (0,1)$ is a non-increasing function.
According to Khinchin's Theorem, 
the inequality \eqref{ineq-Khinchin} has infinitely many solutions
for almost all $x\in\R^d$ if and only if 
$$
{\sum}_{q\ge 1} q^d\psi(q)^d=\infty.
$$
This result raises the problem of estimating the number of solutions for the inequality (\ref{ineq-Khinchin}) satisfying
a specified bound on the denominators, namely analyzing the counting function 
$$
N_T({x}) :=\big|\big\{({p},q) \in \Z^d\times \N:\, 1\le q\le T\;\;\hbox{and}\;\;\hbox{\eqref{ineq-Khinchin} holds } \big\}\big|\,.
$$
The study of this question lead to a major quantitative refinement of Khinchin's Theorem that was established in full generality by W. Schmidt \cite{Sch60} (see also \cite{Le59,Er59} for previous results).
It is natural to embed the set $\Q^d$ as a lattice subgroup in the  space $\mathbb{A}^d$, where $\mathbb{A}$ denotes the rational ad\'eles. Then $N_T({x})$ can be interpreted as the number of lattice points contained in the corresponding domains of $\mathbb{A}^d$. Therefore,
one expects that $N_T({x})$ is approximated by the following volume sum
$$
V_T:=\sum_{1\le q\le T} \vol\big(B\left({x},\psi(q)\right)\big) q^d
=\sum_{1\le q\le T} \big(2\psi(q)\big)^d q^d\,,
$$
where $B(x,\epsilon):=\{{y}\in \R^d:\, \norm{{x} - {y}}_\infty < \epsilon\}.$
Indeed, W.~Schmidt \cite{Sch60} proved that when $V_T\to\infty$, for every $\theta>1/2$,
\begin{equation}\label{eq:schmidt}
N_T(x)= V_T+O_{x,\theta}\left(V_T^{\theta}\right)\quad\quad\hbox{for a.e. $x\in \R^d$}.
\end{equation}


The natural problem of proving analogues of Khinchin's and Schmidt's theorems for rational points on homogeneous algebraic varieties, and in particular on algebraic groups, was raised by S. Lang \cite[p.~189]{L65}.

In previous work \cite{GGN2}, an analogue of Khinchin's Theorem
for rational points on semisimple group varieties was established. The first main goal of the present paper is to establish an analogue of W. Schmidt's asymptotic formula \eqref{eq:schmidt} in this setting. The second main goal is to study systematically the discrepancy of distribution for rational 
points on general semisimple group varieties. 

We turn now to formulate and explain the results in the simplest instance, namely that simply-connected groups defined and almost-simple over $\QQ$.

\subsection{An analogue of W. Schmidt's theorem for group varieties}
Let ${\sf G}\subset \hbox{GL}_n$ be a linear algebraic group defined over $\Q$. For a set of primes $S$, we denote by $\Z_S$ the ring of rational numbers which are integral for every $p\notin S$ (also called $S$-integers), namely whose reduced denominator is divisible only by $p^k, k \ge 0, p\in S$. We consider the group ${\sf G}(\Z_S)$ consisting of matrices whose entries belong to the ring $\ZZ_S$.
The compact open ring of $p$-adic integers in $\QQ_p$ will be denoted $\widehat{\Z}_p$, and ${\sf G}(\widehat{\Z}_p)$ denotes the group of $\widehat{\Z}_p$-points in ${\sf G} (\QQ_p)$. 

It is natural to order the rational numbers in $\ZZ_S$ with respect to the height function
\begin{equation}\label{eq:height0}
\hbox{H}_f (r):=\prod_{p\;{\hbox{\tiny - prime}}} \max(1,\|r\|_p)
=\prod_{p\in S} \max(1,\|r\|_p)\quad\quad \hbox{for $r\in {\sf G}(\Z_S)$,}
\end{equation}
where $\|\cdot\|_p$ denotes the $p$-adic norm on the matrix space $\hbox{Mat}_n(\Q_p)$.
For $x\in {\sf G}(\R)$ and a parameter $b>0$, we consider the inequality
\begin{equation}\label{eq:norm}
\|x-r\|_\infty< \hbox{H}_f(r)^{-b},\quad \hbox{ with $r\in {\sf G}(\Z_S)$.}
\end{equation}
In previous works \cite{GGN1,GGN2},
 the existence of solutions of \eqref{eq:norm}
when $\sf G$ is a connected semisimple algebraic group was investigated.
In particular, when $\sf G$ is almost-simple, simply connected,
and isotropic over $S$, namely ${\sf G}(\QQ_p)$ is non-compact for some $p\in S$, the following was proved. There exist explicit positive exponents $b_1(S)>
b_2(S)$ such that for $b< b_1(S)$,
the inequality \eqref{eq:norm} has infinitely many solutions for almost all $x\in {\sf G}(\R)$, and for $b< b_2(S)$,
the inequality \eqref{eq:norm} has infinitely many solutions for all $x\in {\sf G}(\R)$. In \cite{GGN20} an asymptotic formula for the number of solutions of the inequality  \eqref{eq:norm} was established, for all $x\in {\sf G}(\RR)$ and every $b < b_2(S)$.

We now turn to state (an instance of) our first main result, namely 
an asymptotic formula analogous to the classical estimate \eqref{eq:schmidt}. 
It will be more convenient to 
work instead of \eqref{eq:norm} with an equivalent inequality, namely 
a fixed right invariant Riemannian metric $\rho$ on ${\sf G}(\R)$.  
For a parameter $b>0$ consider the inequality
\begin{equation}\label{eq:metric}
\rho(x,r)\le \hbox{H}_f(r)^{-b},\quad \hbox{ with $r\in {\sf G}(\Z_S)$.}
\end{equation}
As will be shown in Lemma \ref{l:norm-equiv} below, in a ball of fixed radius $R$ centered at a point $x$, the distances in \eqref{eq:norm} and \eqref{eq:metric} are comparable up to a multiplicative constants (depending on $R$).

We define 
\begin{equation}\label{NT}
N_T(x):=\big|\{r\in {\sf G}(\Z_S):\, 1\le \h(r)\le T\;\;\hbox{and}\;\; \hbox{\eqref{eq:metric} holds} \}\big|.
\end{equation}
Note that this counting function is analogous to the one appearing in
\eqref{eq:schmidt}, when the gauge function used is given by $\psi_b(q)=\frac{1}{q^b}$. 
 
We will show that as in \eqref{eq:schmidt} this counting function can be approximated by a suitable adelic volume. We denote by $G_S$ the {\it restricted direct product} of the groups ${\sf G}(\Q_p)$, $p\in S$, with respect to
the compact open subgroups ${\sf G}(\widehat{\Z}_p)$.
 
The definition of the height function \eqref{eq:height0}
extends to the group $G_S$ as 
\begin{equation}\label{eq:height1}
\hbox{H}_f (g):=\prod_{p\in S} \max(1,\|g_p\|_p)\quad\quad \hbox{for $g=(g_p)_{p\in S}\in G_S$.}
\end{equation}
The diagonal embedding 
${\sf G}(\Z_S)\hookrightarrow {\sf G}(\R)\times G_S$
realizes ${\sf G}(\Z_S)$ as lattice subgroup in the product ${\sf G}(\R)\times G_S$. We fix Haar measures 
$m_\infty$ and $m_S$ on the groups ${\sf G}(\R)$ and $G_S$ respectively such that the subgroup ${\sf G}(\Z_S)$ has covolume one in ${\sf G}(\R)\times G_S$ with respect to $m_\infty\times m_S$. 
We consider the volume sum:
$$
V_T:=\sum_{1\le h\le T} m_\infty\big(B(x,h^{-b})\big)m_S\big(\Sigma_S(h)\big),
$$
where 
$$
B(x,\epsilon):=\{y\in {\sf G}(\R): \rho(y,x)< \epsilon\}\quad \hbox{and} \quad\Sigma_S(h):=\{g\in G_S:\, \h(g)=h\}.
$$
We note that because of invariance of the distance $\rho$ and Haar measure $m_\infty$, 
this sum is independent of $x$. The sets $\Sigma_S(h)$ are the {\it height spheres} on the group $G_S$, and they constitute compact open subsets of $G_S$ (and can be empty, a possibility that can certainly occur  for certain values of $h$). 

With this notation, we prove 
the following analogue of W.~Schmidt's result \eqref{eq:schmidt}:

\begin{thm}\label{th:schmidt0}
Let $\sf G$ be a connected simply connected $\Q$-almost simple linear algebraic group defined over $\Q$, and $S$ a finite set of primes
such that $\sf G$ is isotropic for all $p\in S$. Then 
there exists explicit $b_0=b_0(S)>0$ such that
for every parameter $b\in (0,b_0)$,
$$
\big\|N_T-V_T\big\|_{L^2(Q)}\ll_{S,Q} V_T^\theta
$$
with explicit $\theta=\theta(S,b)\in (0,1)$ and an arbitrary bounded measurable subset $Q$ of ${\sf G}(\R)$.
Moreover, for every $\theta'\in (\theta,1)$,
$$
N_T(x)=V_T+O_{S, x,\theta'}\big(V_T^{\theta'}\big)\quad\quad\hbox{for almost every $x\in {\sf G}(\R)$.}
$$ 
\end{thm}

We note that in Theorem \ref{th:schmidt0} we require the assumption that 
$\sf G$ is isotropic (namely ${\sf G}(\QQ_p)$ is non-compact) for every $p\in S$. Our other results below 
will be proved under the weaker condition that $\sf G$ is 
is isotropic for at least one $p\in S$.

Theorem \ref{th:schmidt0} will be proved in Section \ref{sec:schmidt} (see Theorem \ref{ae-disc} and Corollary \ref{cor_NT}). 
We will see that the parameters $b_0$ and $\theta$ can be estimated explicitly in terms of the integrability exponents (see \eqref{eq:q_v} and \eqref{eq:p_v} below) of the relevant automorphic representations.
Moreover, when the automorphic representations 
are known to be tempered and $\sf G$ is unramified over $S$, Theorem \ref{th:schmidt0}
holds with $b_0$ being the divergence exponent of the sum $V_T$ and $\theta=1/2+\eta$
for every $\eta>0$ (see Corollary~\ref{rem:temp}). Hence, in this situation,
we obtain the estimate: for every $\eta>0$,
$$
N_T(x)=V_T+O_{S, x,\eta}\big(V_T^{1/2+\eta}\big)\quad\quad\hbox{for a.e. $x\in {\sf G}(\R)$.}
$$ 
We note that in this case, the exponent of the error estimate of $N_T$ is the square root of the main term. Thus it is of the same quality as Schmidt's theorem stated in \ref{eq:schmidt} above, and can be expected to be best possible.  We refer to Section \ref{sec:best} for further discussion.

\subsection{Strong approximation and discrepancy bounds}
Let ${\sf G}\subset \hbox{GL}_n$ be a connected $\Q$-almost-simple linear algebraic group defined over $\Q$.  In order to render  our introduction more transparent, in the present section we assume that ${\sf G}$ is simply connected, and defer for later a discussion of the general case. Then it is known
that $\sf G$ satisfies the Strong Approximation Property \cite[\S7.4]{PlaRa}. This property implies that for every set $S$ of primes such that ${\sf G}$ is isotropic over $\Q_p$
for some $p\in S$, the diagonal embedding 
$$
{\sf G}(\Z_S)\hookrightarrow {\sf G}(\R)\times \cI^S,\quad\quad\hbox{where $\cI^S:=\prod_{p\notin S} {\sf G}(\widehat{\Z}_p)$,}
$$
is dense. 

For a measurable subset $E$ of ${\sf G}(\R)$ and $x\in {\sf G}(\R)$, we set 
$$
E(x):=Ex.
$$
We will analyze the distribution of
the rational points ${\sf G}(\Z_S)$ in the subsets $E(x)\times W$,
 where $W$ are compact open subsets of $\cI^S$.
We fix the invariant probability measure $m^S$ on $\cI^S$, a Haar measure $m_S$ on $G_S$,
and a Haar measure $m_\infty$ on ${\sf G}(\R)$.
We recall that under the diagonal embedding,
${\sf G}(\Z_S)$ is a discrete subgroup of finite covolume in 
${\sf G}(\R)\times G_S$. We normalize the measure $m_\infty$ and $m_S$ so that the subgroup ${\sf G}(\Z_S)$ has covolume one in the product. We set
\begin{align*}
R_S(h) &:= \left\{\gamma\in {\sf G}(\Z_S): \, \hbox{H}_f(\gamma)\le h\right\},\\
B_S(h)&:=\left\{g\in G_S: \, \hbox{H}_f(g)\le h\right\}\quad\hbox{and}\quad 
v_S(h):=m_S\big(B_S(h)\big).
\end{align*}
We will be interested in estimating the cardinality
$\big|R_S(h)\cap (E(x)\times W)\big|$. However, this cardinality
might be infinite (for instance, when $E$ contains a coset of 
the group ${\sf G}(\Z)$). To address this issue, we define 
$$
\mathcal{N}(E):=\big|\{\gamma\in {\sf G}(\mathbb{Z}):\, m_\infty(E\cap \gamma E) > 0\}\big|,
$$
and assume that $\mathcal{N}(E)<\infty$. For example, if $E$ is bounded, this is always the case.

We define the discrepancy of the  rational points ${\sf G}(\Z_S)$ as
\begin{equation}\label{disc-def}
\mathcal{D}\big(R_S(h), E(x)\times W\big):=\left| \frac{|R_S(h)\cap (E(x)\times W)|}{v_S(h)} - m_\infty(E)m^S(W) \right|.
\end{equation}  
Remarkably, we show that the discrepancy can be estimated for {\it general measurable
domains} $E$ of finite measure, as follows. 

\begin{thm}\label{th:intr1}
Let $\sf G$ be a connected simply connected $\Q$-almost-simple linear algebraic group defined over $\Q$, and let $S$ be a set of primes
such that $\sf G$ is isotropic for some $p\in S$.	
Then there exists $\mathfrak{k}_S({\sf G})>0$ such that for every measurable subset $E$ of  ${\sf G}(\R)$ with finite measure satisfying $\mathcal{N}(E)<\infty$, for 
every compact open subset $W$ of $\cI^S$,  and for every $\eta > 0$.
$$
\Big\|\mathcal{D}\big(R_S(h), E(\cdot)\times W\big)\Big\|_{L^2(Q)}\ll_{S,Q,\eta} \mathcal{N}(E)^{1/2} m_\infty(E)^{1/2}\,m^S(W)^{1/2}\, v_S(h)^{-\mathfrak{k}_S({\sf G})+\eta}
$$
for every bounded measurable subset $Q$ of ${\sf G}(\R)$.

\end{thm}

We note that the spectral exponent $\mathfrak{k}_S({\sf G})$ will be given an explicit form in our discussion of the proof of Theorem \ref{th:intr1} in \S 4.

Using the $L^2$-bound established in Theorem \ref{th:intr1},
we also deduce an almost sure estimate on the discrepancy:

\begin{thm}\label{cor:aa}
With notation as in Theorem \ref{th:intr1},	for every 
$0 < \mathfrak{k} < \mathfrak{k}_S({\sf G}) $, and for almost every $x\in {\sf G}(\R)$, and for every $\eta > 0$ 
$$
\mathcal{D}\big(R_S(h), E(x)\times W\big) \ll_{S,E,W,x,\mathfrak{k},\eta}\, \big(\log v_S(h)\big)^{3/2+\eta} v_S(h)^{-\mathfrak{k}}.
$$
\end{thm}

\begin{rem}\label{rm:meas-zero}
{\rm 
We note that when $m_\infty(E)=0$, Theorem \ref{cor:aa} provides the almost sure upper bound 
$$|R_S(h)\cap (E(x)\times W)| \ll_{S,E,W,x,\mathfrak{k},\eta}\big(\log v_S(h)\big)^{3/2+\eta} v_S(h)^{1-\mathfrak{k}}$$
for the number of rational points of height bounded by $h$ in $E$ satisfying the congruence constraints defined by $W$. Therefore, it amounts to a general non-concentration phenomenon for rational points. As an example, $E\subset G_\infty $ can be any smoothly embedded submanifold of positive co-dimension satisfying $\mathcal{N}(E)< \infty$.
}
 \end{rem}
\vspace{0.2cm}

Clearly, one has to impose additional assumptions on the sets $E$
to expect an estimate for the discrepancy $\mathcal{D}\big(R_S(h), E(x)\times W\big)$
that is valid for all $x$. 
Indeed, when the sets $E$ satisfy a suitable regularity property, 
we establish such a pointwise bound, as follows.

We say that a measurable subset $E$ of 
${\sf G}(\R)$ is \emph{right-stable} if  
\begin{equation}\label{eq:well}\tag{RS}
m_\infty(E_\epsilon^+\backslash E_\epsilon^-)\ll \epsilon \quad\hbox{ for every $\epsilon\in (0,\epsilon_0)$,} 
\end{equation}
where 
$$
E_\epsilon^+:=E B(e,\epsilon) \quad\hbox{and}\quad 
E_\epsilon^-:=\{x\in {\sf G}(\R):\, x B(e,\epsilon)\subset E \}.
$$
\begin{rem}\label{rm:right-stable}
{\rm 
The definition of right-stability is motivated by the notion of well-roundedness from \cite{GN12}, but is considerably more general. Note that many sets of measure zero (including finite sets) are right-stable.  These include, for example,  the intersection of smoothly embedded positive-codimenion submanifolds of $G$ with a norm ball. 
}
\end{rem}

Our pointwise error estimate will now depend on the dimension $d:=\dim_\R({\sf G}(\R))$.

\begin{thm}\label{th:wr}
Let $\sf G$ be a connected simply connected $\Q$-almost-simple linear algebraic group defined over $\Q$, and $S$ a set of primes
such that $\sf G$ is isotropic for some $p\in S$.	
Let $E$ be a right-stable finite-measure subset of ${\sf G}(\R)$ satisfying $\mathcal{N}(E_{\epsilon_0}^+)<\infty$, and $W$ compact open subset of $\cI^S$. Then for $x\in {\sf G}(\R)$ for for every 
$0 < \mathfrak{k}= \mathfrak{k}_S({\sf G})-\eta < \mathfrak{k}_S({\sf G})$ 
\begin{align*}
\big|R_S(h)\cap (E(x)\times W)\big|=&\,m_\infty(E)m^S(W)v_S(h)\\
&+ O_{S,E,x,\eta}\Big( m^S(W)^{(d+1)/(d+2)}\, v_S(h)^{1-2\mathfrak{k}/(d+2)}\Big)\,,
\end{align*}
provided that 
$m^S(W)\gg_\eta \, v_S(h)^{-2\mathfrak{k}}.$
This condition is equivalent  to the error term in the estimate being bounded by the main term. 

Explicity, if the volume growth satisfies $v_S(h)\gg_{S,\eta^\prime} h^{\mathfrak{a}}$, with $0< \mathfrak{a}=\mathfrak{a}_S({\sf G})-\eta^\prime$, then the estimate holds provided the height $h$ satisfied $h\gg_{\eta^\prime, \eta}  (1/m^S(W))^{1/2\mathfrak{a}\mathfrak{k}} =h_{\eta,\eta^\prime}(S,W)$.
Moreover, the above  estimate is uniform for $x$ ranging in compact subsets of ${\sf G}(\R)$.
\end{thm}

In particular, when $m_\infty(E) > 0$, it follows that for such $h$,
$R_S(h)\cap (E(x)\times W)\ne \emptyset$, 
 namely the domain $E(x)\times W$ contains a point in ${\sf G}(\Z_S)$ (namely an $S$-integral point) with height at most $h$. 
 
 We note that Theorem \ref{th:wr} and its more general version stated in \S 4 generalize several earlier results, including \cite{Cl02},  \cite{Du03}, \cite{Oh05} and \cite[Thm. 1.8]{BO}. 

The method of the proof of Theorem \ref{th:wr}
can be used to establish bounds on discrepancy which are uniform over 
variable families of sets. To demonstrate the utility of this fact, we analyze the  
discrepancy with respect to the entire family of Riemannian balls $B(x,l)$ in ${\sf G}(\R)$, with $x\in G$ and $0< l < l_0$,  
and establish explicit discrepancy estimates for all points in the group, at arbitrary small scales.

\begin{thm}\label{th:balls}
Let $\sf G$ be a connected simply connected $\Q$-almost-simple linear algebraic group defined over $\Q$, and $S$ a set of primes
such that $\sf G$ is isotropic for some $p\in S$. Fix a compact open subset $W$ of $\cI^S$, and $0 < \mathfrak{k}= \mathfrak{k}_S({\sf G})-\eta < \mathfrak{k}_S({\sf G})$. For a suitable $\ell_0 >0$, we set for $0 < \ell < \ell_0$,   
$$
\mathcal{E}_{\ell,W}(h):= m_\infty(B(e,\ell))^{d/d+2}m^S(W)^{(d+1)/(d+2)}v_S(h)^{1-2\mathfrak{k}/(d+2)}\,,
$$
Then for every $x\in {\sf G}(\R)$
$$
\big|R_S(h)\cap (B(x,\ell)\times W)\big|= m_\infty(B(e,\ell))m^S(W)v_S(h)+O_{S,x, \eta}\Big (\mathcal{E}_{\ell,W}(h)\Big),
$$
provided that $m_\infty(B(e,\ell))^2m^S(W)\gg_\eta v_S(h)^{-2\mathfrak{k}}$. We note that this condition is equivalent to the error term in the estimate being bounded by the main term. 
%

Explicity, if the volume growth satisfies $v_S(h)\gg_{S,\eta^\prime} h^{\mathfrak{a}}$, with $0< \mathfrak{a}=\mathfrak{a}_S({\sf G})-\eta^\prime$, then the estimate holds provided the height $h$ satisfies
$$h\gg_{S,\eta,\eta^\prime}\ell^{-d/\mathfrak{k}\mathfrak{a}} m^S(W)^{-(d+2)/2\mathfrak{k}\mathfrak{a}}\,.$$ 
Moreover, this estimate is uniform for $x$ ranging in compact subsets of ${\sf G}(\R)$.
\end{thm}

The results stated in this section will be proven in Section \ref{sec:dis_simp}.

		We note that the effective estimates on discrepancy stated in Theorem \ref{th:wr} and Theorem \ref{th:balls} can also be viewed as establishing an effective count for the number of rational solutions of intrinsic Diophantine inequalities. This latter problem can be reduced to a lattice counting problem, and was given a short and simple solution in \cite{GGN20}, under more restrictive hypotheses than those of Theorem \ref{th:wr} and Theorem \ref{th:balls} .  Indeed ${\sf G}(\mathbb{Z}_S)$ 
		can be viewed as a lattice subgroup in ${\sf G}(\mathbb{R})\times G_S$, and when  
		a (variable) family of balls $D(x, \vare)\times B_S(h)$ is suitably well-rounded (as a function of $\vare$), the effective solution of the lattice point counting problem in  \cite{GN12} can be applied.
		
		Theorem \ref{th:wr} and Theorem \ref{th:balls} generalize the results in \cite{GGN20} in several respects, as follows.  
		\begin{itemize}
		\item The set $E$ need only be a measurable set of positive finite measure satisfying right-stability, a considerably weaker condition than well-roundedness,   
\item  A stronger error bound for the discrepancy $\mathcal{D}\big(R_S(h), B(x,l)\times \cI^S\big) $ is established. 
		For comparison, in the case of Riemannian balls $B(x,\ell)$, the method of the present paper gives  the bound, for any $\eta > 0$ 
		$$
\ll_{S,x,\eta} m_\infty(B(e,\ell))^{d/(d+2)}
v_S(h)^{(-2\mathfrak{k}/(d+2))+\eta}.
$$
	while the bound established in \cite[Thm. 1.3]{GGN20} gives, for a set $S$ of unramified  primes (see the discussion following Theorem \ref{th:mean_sc2} below), and for the choice $W=\cI^S$ (namely in the absence of congruence conditions)
		$$
		\ll_{S,x,\eta} m_\infty(B(e,\ell))^{d/(d+1)}\, v_S(h)^{(-\mathfrak{k}/(d+1))+\eta}.
		$$	
		\item  An arbitrary congruence constraint is allowed on the rational points in ${\sf G}(\mathbb{Z}_S)$ involved in the approximation process, given by an arbitrary compact open subset $W$ of $\cI^S$. 
	\end{itemize}
	Furthermore, we will generalize Theorems
\ref{th:intr1}--\ref{th:balls} in two additional important respects, namely  we will consider {\it every} connected almost $\QQ$-simple algebraic group, not only simply-connected one, and we will consider {\it every} non-empty subset $S\subset P$ of prime over which the $\QQ$-group ${\sf G}$ is isotropic, including subsets which contain ramified primes.  Each of these extensions  require elaborate arguments using structure theory of ad\'ele groups and spectral results in their automorphic representations. Let us now turn to formulate the corresponding results.

\subsection{Discrepancy bounds for general $\QQ$-groups}

We now allow the group $\sf G$ to be a general $\QQ$-almost simple group, not necessarily simply connected. In this case, the strong approximation property fails, and in particular the embedding of ${\sf G}(\Z_S)$ in ${\sf G}(\R)$ is not dense, typically. Nonetheless, one can show (cf. Corollary \ref{c:dense}) that the closure of 
${\sf G}(\Z_S)$ in  ${\sf G}(\R)$ is a finite index subgroup of ${\sf G}(\R)$.
Since $\overline{{\sf G}(\Z_S)}$ decomposes as a finite union of cosets of 
${\sf G}(\R)^0$, the connected component of identity in ${\sf G}(\R)$,
it is sufficient to analyze the distribution of ${\sf G}(\Z_S)$-points in ${\sf G}(\R)^0$.

To state our results precisely, we need to take into account the contribution of  automorphic characters, as follows. Let $\mathcal{X}({\sf G}, \cI_f)$ denote the set of continuous unitary characters $\chi$ on the ad\'ele group ${\sf G}(\mathbb{A}_\Q)$ such that 
$\chi({\sf G}(\Q))=1$ and $\chi({\sf G}(\widehat{\Z}_p))=1$ for all $p$.
We denote by $G^{\ker}$ the joint kernel of this (finite) set of characters of ${\sf G}(\mathbb{A}_\Q)$, and note that $G^{\ker}$ is a finite index subgroup of ${\sf G}(\mathbb{A}_\Q)$
(see, for instance, \cite[Lemma~4.4]{GGN1}).

We set 
$$
G_S^{\ker}:=G_S\cap G^{\ker}\quad \hbox{and}\quad
v_S^{\ker}(h):= m_S \Big( \big\{g\in G^{\ker}_S:\, \h(g)\le h   \big\} \Big).
$$
Let $m_\infty^0$ be the Haar measure on ${\sf G}(\R)^0$ normalized so that 
the intersection of ${\sf G}(\Z_S)$ with ${\sf G}(\R)^0\times G_S^{\ker}$
has covolume one in ${\sf G}(\R)^0\times G_S^{\ker}$ with respect to the measure $m_\infty^0\times m_S$. For $E\subset {\sf G}(\R)^0$, we consider the discrepancy  

\begin{equation}\label{disc-non-SC}
\mathcal{D}\big(R_S(h), E\big):=\left| \frac{|R_S(h)\cap E|}{v^{\ker}_S(h)} - m^0_\infty(E)\right|.
\end{equation}

Note that the in the foregoing expression we did not impose any congruence conditions on the rational points involved, as we did in the simply-connected case,  although our methods certainly allow for this possibility. But since in the non-simply-connected case certain congruence obstructions typically arise, a complete analysis of the discrepancy of rational points subject to congruence conditions requires considerably more notation and further discussion, which we will avoid in order to keep the exposition more accessible.  

We shall show that with this setup an analogue of Theorem \ref{th:intr1} holds:

\begin{thm}\label{th:intr1_2}
	Let $\sf G$ be a connected $\Q$-almost-simple linear algebraic group defined over $\Q$, and $S$ a set of primes
	such that $\sf G$ is isotropic for some $p\in S$.	
	Then there exists $\mathfrak{k}_S({\sf G})>0$ such that for every measurable subset $E$ of ${\sf G}(\R)^0$ of finite measure satisfying $\mathcal{N}(E)<\infty$, and any $\eta > 0$ 
	$$
	\Big\|\mathcal{D}\big(R_S(h), E(\cdot)\big)\Big\|_{L^2(Q)}\ll_{S,Q,\eta} \mathcal{N}(E)^{1/2} m^0_\infty(E)^{1/2}\,m^S(W)^{1/2}\, v^{\ker}_S(h)^{-\mathfrak{k}_S({\sf G})+\eta}
	$$
	for every bounded measurable subset $Q$ of ${\sf G}(\R)^0$.
\end{thm}
Theorem \ref{th:intr1_2} will be proved in \S  \ref{sec:gen}.
\begin{rem}\label{subtle}
{\rm 
The estimate of Theorem \ref{th:intr1_2} is similar in form to
Theorem \ref{th:intr1}, but there is an important, if subtle, difference between them. The normalization of Haar measure $m_\infty^0$ 
is typically different than that of  $m_\infty$, and the volume functions 
$v_S(h)$ and $v_S^{\ker}(h)$, while comparable, are typically also different. 
Therefore the existence of automorphic characters influences the size of the main term in the asymptotic formulas associated with (\ref{disc-def}) and (\ref{disc-non-SC}), replacing  
$v_S(h)m_\infty(E)$ which arise in the simply connected case by $v_S^{\ker}(h)m_\infty^0(E)$. 
}
\end{rem}

\begin{rem}
{\rm
We note that the closure $\overline{{\sf G}(\Z_S)}$ in ${\sf G}(\R)$
is a union of finitely many cosets $\gamma_i {\sf G}(\R)^0$ with $\gamma_i\in {\sf G}(\Z_S)$. Given a subset $E$ of $\overline{{\sf G}(\Z_S)}$, we may decompose it as 
$E=\sqcup_i E_i$ with  $E_i:=E\cap \gamma_i {\sf G}(\R)^0$ and estimate the discrepancy separately for each of the subset $E_i$. 
Hence, our method can be used more generally to analyze discrepancy for 
subsets  $E\subset \overline{{\sf G}(\Z_S)}$.
}	
\end{rem}

Once the mean-square bound for $\mathcal{D}\big(R_S(h), E(x))$ 
has been established, we can also prove generalizations of Theorems
\ref{th:intr1}--\ref{th:balls} with obvious modifications, which will be explained further in  \S  \ref{sec:gen}.

In summary, in the present section we gave an account of our results for almost-simple linear algebraic groups defined over $\QQ$, to simplify the presentation. Our methods, however, are completely general, and from now on we will turn to develop discrepancy estimates when
\begin{itemize}
\item the ground field $K$ is arbitrary algebraic number field, namely a finite dimensional extension of $\QQ$, 
\item the group $\sf G$ is an arbitrary connected semisimple linear algebraic group defined and simple ver $K$, not necessarily simply-connected, 
 \item the set of places $S$ of $K$ may contain ramified places, 
 \item the approximating rational elements are subject to an arbitrary congruence constraint when the group is simply connected. 
\end{itemize}

\subsection*{Acknowledgements}
A.G. was supported by SNF grant 200021--182089, and 
A.N. was supported by ISF Moked Grant 2919-19.

\section{Notation and preliminary results}\label{sec:not}

\subsection{Algebraic groups over number fields}\label{sec:alg-gps}
We start by reviewing basic properties of semisimple algebraic groups over number fields (cf. \cite[Ch.~3--5]{PlaRa}).
Throughout the paper, $K$ denotes an algebraic number field (namely, a finite extension of the field of rationals $\QQ$). Let $V_K$ be the set of normalized absolute values $|\cdot|_v$ of the field $K$. 
The set of absolute values  decomposes as 
$$
V_K=V_K^\infty\sqcup V_K^f,
$$
where $V_K^\infty$ is the finite set of Archimedian absolute values and 
$V_K^f$ the set of non-Archimedian absolute values.
For $v\in V_K$,
we write $K_v$ for the corresponding completion of $K$, and
when $v\in V_K^f$, we denote by 
$$
O_v:=\{x\in K_v:\, |x|_v\le 1\}
$$ 
the ring of integers in $K_v$.
For $S\subset V^K_f$, we write
$$
O_S:=\{x\in K:\, |x|_v\le 1\hbox{ for $v\notin S$} \}\footnote{Note that $O_{\{v\}}\neq O_v$, the ring of integers  just defined !} \,.
$$
for the ring of $S$-integers, namely elements in $K$ that are integral w.r.t. every completion except possibly those in $S$. In particular, $O=O_\emptyset$ denotes the ring of integers in $K$.

We denote by 
$$
{\mathbb A}_K:=\big\{(x_v)_{v\in V_K}:\, |x_v|_v\le 1\hbox{ for almost all $v$}\big\}
$$ the ring of ad\'eles of $K$, which is
the restricted direct product of $K_v$, $v\in V_K$, with respect to compact open subrings $O_v\subset K_v$, $v\in V_K^f$.

Let ${\sf G}\subset \hbox{SL}_n$ be 
a connected semisimple linear algebraic group  defined over the field $K$. 
The focus of our investigation is the distribution 
of the set of $S$-integral points
$$
\Gamma_S:={\sf G}(O_S).
$$
We denote by 
$$
G_v:={\sf G}(K_v)\quad \hbox{for $v\in V_K,$}
$$
the locally compact groups of $K_v$-points of $\sf G$ equipped with the topology defined by the corresponding absolute values $|\cdot|_v$. 
We say that $\sf G$ is \emph{isotropic} over $v$ if $G_v$ is non-compact,
and $\sf G$ is \emph{isotropic} over $S\subset V_K$ if $G_v$ is non-compact
for at least one $v\in S$. 

For $v\in V_K^f$, we consider
$$
{\sf G}(O_v):=\{g\in G_v:\, \|g\|_v\le 1\},
$$
which is a compact open subgroup of $G_v$. 
We introduce compact groups
$$
\cI_S:={\prod}_{v\in S}\; {\sf G}(O_v)\quad\quad \hbox{and} \quad\quad \cI^S:={\prod}_{v\in V_K^f\backslash S}\; {\sf G}(O_v). 
$$
For $S\subset V_K$,
we write 
$$
G_S:={\prod}_{v\in S}'\, G_v=\big\{(g_v)_{v\in S}:\, g_v\in G_v,\, \|g_v\|_v\le 1\hbox{ for almost all $v$}\big\}.
$$
for the restricted direct product of the groups $G_v$ with
respect to the compact open subgroups ${\sf G}(O_v)$, $v\in S\cap V_K^f$.
Then $G_S$ is a locally compact group. For instance, ${\sf G}(\mathbb{A}_K)=G_{V_K}$ is the ad\'ele group associated to ${\sf G}$. To simplify notation, we also write 
$$
G_\infty:=G_{V_K^\infty}\quad\hbox{and}\quad G_f:=G_{V_K^f}.
$$
We recall that when $\sf G$ is simply connected, the group $G_\infty$ is connected
with respect to the Euclidean topology
(cf. \cite[Ch.~7, Prop.~7.2]{PlaRa}). In general, $G_\infty$ has finitely many connected components, and we denote by $G_\infty^0$
the connected component of the identity in $G_\infty$.

For $v\in V_K$, we denote by $m_{v}$ the Haar measure on $G_v$, and when $v\in V_K^f$ we normalize it so that $m_v({\sf G}(O_v))=1$. Then for a subset $S\subset V_K^f$, the product measure $m_S:=\prod_{v\in S} m_v$ defines a Haar measure on $G_S$ such that $m_S(\cI_S)=1$.
We also denote by $m_\infty$ a Haar measure on $G_\infty$.
Under the diagonal embedding 
$\Gamma_S\hookrightarrow G_\infty\times G_S$,
the group $\Gamma_S$ is a discrete subgroup 
with finite covolume in $G_\infty\times G_S$ (cf. \cite[Ch.~5]{PlaRa}).
We normalize the measure $m_\infty$, so that $\Gamma_S$ has covolume one with respect $m_\infty\times m_S$.


We recall the Strong Approximation Property \cite[\S 7.4]{PlaRa}:

\begin{thm}[Strong Approximation]\label{th:strong_app}
	Let $\sf G$ be a simply connected $K$-simple algebraic group defined over $K$.  Then if $\sf G$ is isotropic over $S\subset V_K^f$,
	the image of $\Gamma_S$ with respect to the embedding 
	$\Gamma_S\hookrightarrow G_\infty\times \cI^S$
	is dense. More generally, the embedding 
	${\sf G}(K)\hookrightarrow G_\infty\times G_{V_K^f\backslash S}$ has dense image.
\end{thm} 

This result fails if the group $\sf G$ is not simply connected even for the  
embedding $\Gamma_S\hookrightarrow G_\infty$. Nonetheless, one can deduce the following result about the closure: 

\begin{cor}\label{c:dense}
	Let $\sf G$ be a connected $K$-simple algebraic group defined over $K$.	Then if $\sf G$ is isotropic over $S\subset V_K^f$, 
	the closure of ${\Gamma_S}$ in $G_\infty$ 
	is a finite index open subgroup. In particular, $\overline{\Gamma_S}\supset G_\infty^0$.
\end{cor}

\begin{proof}
	Let us consider the simply connected cover ${\sf p}:\tilde {\sf G}\to {\sf G}$.
	Then the groups :
	$$
	U_v:={\sf p}^{-1}({\sf G}(O_v))\cap \tilde {\sf G}(O_v) \quad \hbox{ with $v\in V_K^f$}\,,
	$$
	are compact open subgroups of $\tilde G_v$. It follows from Theorem \ref{th:strong_app}  that the image of the group
	$$
	\Gamma:=\tilde {\sf G}(O_S)\cap \left(\tilde G_\infty\times{\prod}_{v\in V_K^f\backslash S} U_v\right)
	$$
	is dense in $\tilde G_\infty$. Since ${\sf p}(\Gamma)\subset \Gamma_S$,
	it follows that the closure $\overline{\Gamma_S}$ in $G_\infty$ contains
	${\sf p}(\tilde G_\infty)$ which is an open and closed subgroup of finite index in $G_\infty$ (cf. \cite[\S 3.2]{PlaRa}).  
\end{proof}

\subsection{Automorphic representations} \label{sec:auto}
Let $\sf G$ be a connected $K$-simple algebraic group defined
over a number field $K$. We consider the Hilbert space 
$$
\mathcal{H}_{\sf G}:=L^2({\sf  G}(\mathbb{A}_K)/{\sf G}(K))
$$
consisting of square-integrable functions on the space 
${\sf  G}(\mathbb{A}_K)/{\sf G}(K)$ equipped with the invariant probability measure $\mu$.
Let
$$
\mathcal{H}^0_{\sf G}:=L_0^2\big({\sf  G}(\mathbb{A}_K)/{\sf G}(K)\big)=\left\lbrace
\phi\in \mathcal{H}_{\sf G}:\, \int_{{\sf  G}(\mathbb{A}_K)/{\sf G}(K)} \phi\, d\mu=0 
 \right\rbrace .
$$
A continuous unitary  character $\chi$ of ${\sf  G}(\mathbb{A}_K)$ is called {\it automorphic} if 
$\chi({\sf G}(K))=1$. Then $\chi$ can be considered as an element of $\mathcal{H}_{\sf G}$.
We denote by $\mathcal{H}_{\sf G}^{00}$ the subspace of $\mathcal{H}_{\sf G}$ orthogonal to all automorphic
characters. We note that when $\sf G$ is simply connected 
there no non-trivial authomorphic characters and $\mathcal{H}_{\sf G}^{00}=\mathcal{H}_{\sf G}^{0}$.

We fix a choice of maximal compact subgroups $U_v$ of $G_v$
that coincide with ${\sf G}(O_v)$ for all but finitely many $v\in V_K^f$.
Furthermore, for all but finitely many $v$, these subgroups satisfy:
\begin{itemize}
	\item[(i)] $U_v$ is a hyperspecial, good maximal compact subgroup of $G_v$,
	\item[(ii)] the group $\sf G$ is unramified over $K_v$ (that is,
	$\sf G$ is quasi-split over $K_v$ and split over an unramified extension of $K_v$).
\end{itemize}
For the remaining finite set of places $v$, we fix a good special maximal compact subgroup $U_v$ of $G_v$. For any  subset $S\subset V^f_K$, we set
$$
U_S:={\prod}_{v\in S}\; U_v\quad \hbox{and} \quad U^S:={\prod}_{v\in V_K^f\backslash S}\; U_v. 
$$
If for all $v\in S$, we have  $U_v ={\sf G}(O_v)$, then 
we say that $G$ is \emph{unramified} over $S$.

We denote by $\pi_v^{\rm aut}=\pi_v$ the unitary representation of the group $G_v$
on the space $\mathcal{H}_{\sf G}$. 

For places $v\in V_K^f$,
the \emph{spherical integrability exponent} of the representations
$\pi_v$, with respect to the subgroup $U_v$, is defined by
\begin{equation}\label{eq:q_v}
\mathfrak{q}_v({\sf G}):=\inf\left\{ q \ge 2:\, \begin{tabular}{l}
$\forall$\hbox{ $U_v$-inv. $\phi\in \mathcal{H}^{00}_{\sf G}$}\\
\hbox{$\left< \pi_v(g)\phi,\phi\right>\in L^q(G_v)$}
\end{tabular}
\right\}.
\end{equation}
It is a fundamental result in the theory of automorphic representations 
that the integrability exponents $\mathfrak{q}_v({\rm G})$ are finite,
and moreover, $\mathfrak{q}_v({\rm G})$ is uniformly bounded over $v$, see \cite{Cl03}.
These exponents can be estimated in terms of 
the Satake parameters of the corresponding spherical automorphic
representations. We refer to \cite{S05},\cite{C07} for surveys of some of these  results.

More generally, for $S\subset V_K^f$, 
we denote by $\pi_S^{\rm aut}=\pi_S$ the unitary representation 
of the group $G_S$ on the space $\mathcal{H}_{\sf G}$.
We define 
\begin{equation}\label{eq:q_v}
\mathfrak{q}_S({\sf G}):=\inf\left\{ q \ge 2:\, \begin{tabular}{l}
$\forall$\hbox{ $U_S$-inv. $\phi\in \mathcal{H}^{00}_{\sf G}$}\\
\hbox{$\left< \pi_S(g)\phi,\phi\right>\in L^q(G_S)$}
\end{tabular}
\right\}.
\end{equation}
  This integrability exponent is also finite  and it can be estimated in terms of 
the exponents $\mathfrak{q}_v({\sf G})$, $v\in S$, see \cite[Cor.~3.5]{GGN1}.

Given a strongly continuous unitary representation $\pi:G_S\to \mathcal{U}(\mathcal{H})$ 
on a Hilbert space $\mathcal{H}$ and  a finite Borel measure $\nu$ on $G_S$, 
we define the averaging operator
$$
\pi_S(\nu):\mathcal{H}\to \mathcal{H}:\, \phi\mapsto \int_{G_S} \pi(g)\phi\, d\nu(g).
$$
We recall the following estimate on the norm of the averaging operators:

\begin{thm}[\cite{GGN1}, Prop. 3.8]\label{th:norm}
	Let $\beta$ be a Haar-uniform probability measure supported on a $U_S$-bi-invariant
	bounded subset $B$ of $G_S$. 
	Then 
	$$
	\left\|\pi_S^{\rm aut}(\beta)|_{\mathcal{H}_{\sf G}^{00}}\right\|\ll_{S,\eta} m_{S}(B)^{-\frac{1}{\mathfrak{q}_S({\rm G})}+\eta}\quad\hbox{ for all $\eta>0$.}
	$$	
\end{thm}
We note that although the measure $m_S$ in \cite{GGN1} was normalized differently 
(namely, so that $m_S(U_S)=1$), this gives the same bound up to a multiplicative constant.

\vspace{0.2cm}

An important fact,  underlying our considerations below regarding ramified places, is that this result also holds for more general averaging operators. The crucial ingredient here
is finiteness of a more general integrability exponent, which we define as follows :
\begin{equation}\label{eq:p_v}
\mathfrak{p}_S({\sf G}):= \inf \left\{p\ge 2:\, 
\begin{tabular}{c}
$\forall$\hbox{ $\phi_1,\phi_2$ in a dense subspace of $\mathcal{H}_{\sf G}^{00}$}\\
\hbox{$\left< \pi_S^{\rm aut}(g)\phi_1,\phi_2\right>\in L^p(G_S)$}
\end{tabular}
\right\}.
\end{equation}
One says that a unitary representation $\pi:G_S\to \mathcal{U}(\mathcal{H})$ is {\it $L^p$-integrable}, if for $\phi_1,\phi_2$ in a dense subset of $\mathcal{H}$,
the functions $g\mapsto \left< \pi(g)\phi_1,\phi_2\right>$ is in $L^p(G_S)$. Thus,
$$
\mathfrak{p}_S({\sf G})= \inf \left\{p\ge 2:\, 
\hbox{ $\pi_S^{\rm aut}|_{\mathcal{H}_{\sf G}^{00}}$ is $L^p$-integrable}
\right\}.
$$
It was proved  in \cite{GMO} that this exponent is finite provided 
that $\sf G$ is either simply connected or adjoint (see \cite[Th.~3.20 and Th.~3.7]{GMO}). 
We shall show in Theorem \ref{th:integrability} below that the exponent is finite for  general  $K$-simple groups.

We define:
$$
\mathfrak{n}_S({\sf G}):=
\left\{
\begin{tabular}{ll}
\hbox{the least even integer $\ge \mathfrak{p}_S({\sf G})/2$, if $\mathfrak{p}_S({\sf G})>2$,}\\
\hbox{1, if $\mathfrak{p}_S({\sf G})=2$.}
\end{tabular}
\right.
$$
With this notation, we have:

\begin{thm}\label{th:norm2}
	Let $W_S$ be a compact open subgroup of $G_S$, and $\beta$ a Haar-uniform 
	probability measure supported on a $W_S$-bi-invariant bounded subset $B$ of $G_S$.
	Then
	$$
	\left\|\pi_S^{\rm aut}(\beta)|_{\mathcal{H}_{\sf G}^{00}}\right\|\ll_{\eta,W_S} m_S(B)^{-\frac{1}{4\mathfrak{n}_S({\sf G})}+\eta} \quad\hbox{	for all $\eta>0$.}
	$$
\end{thm}

\begin{proof} 
	The proof is a generalization of the proof of \cite[Cor.~6.7]{GN12}, and so we only provide an outline and refer 
	\cite[\S 6]{GN12} for further details.
	The main idea, which originated in \cite{CHH} and  \cite{N98},
	is to observe that a suitable tensor power of $\pi_S$
	restricted to $\mathcal{H}_{\sf G}^{00}$ 
	is weakly contained in (a multiple of) the regular representation $\lambda_S$ of $G_S$ on $L^2(G_S)$, and then use a generalization of the Kunze--Stein convolution
	inequality valid in $L^2(G_S)$. More precisely, 
	the representation $(\pi_S|_{\mathcal{H}_{\sf G}^{00}})^{\otimes\mathfrak{n}_S({\sf G})}$ is 
	$L^p$-integrable for all $p>2$, so that it is weakly contained in
	the regular representation $\lambda_G$ (by \cite{CHH}), which allows to deduce (by \cite{N98}) 
	that for any probability density $\beta^\prime$ on $G_S$ 
\begin{equation}\label{eq:11}
\left\|\pi_S^{\rm aut}(\beta^\prime)|_{\mathcal{H}_{\sf G}^{00}}\right\|\le
\|\lambda_S(\beta^\prime)\|^{\frac{1}{\mathfrak{n}_S({\sf G})}}.
\end{equation}
Now let 
	$$
	 B^\prime:=U_S B U_S.
	$$
	and denote by $ \beta^\prime$ the Haar-uniform probability measure
	supported on $B^\prime$. Since $U_S\cap W_S$ has finite index
	in both $U_S$ and $W_S$, it is clear that there exist $c_1,c_2>0$, depending only on $W_S$, such that
	$$
	c_1\, m_S(B^\prime)\le m_S(B)\ll c_2\, m_S( B^\prime),
	$$
	so that 
	\begin{equation}\label{eq:2}
	\|\lambda_S(\beta)\|\le c\, \|\lambda_S(\beta^\prime)\|
	\end{equation}
	for some $c>0$.
	Finally, using the decomposition of $m_S$ with respect to the Iwasawa 
	decomposition on $G_S$, we deduce (cf. \cite[Thm.~6.6]{GN12}) that
	$\|\lambda_S( \beta^\prime)\|$ can be estimated in terms of the Harish Chandra function on $G_S$.  The Harish Chandra function in this case is $L^{4+\eta}$-integrable for all $\eta>0$ 	(by \cite[Prop.~6.3]{GN12}), and we deduce that
	\begin{equation}\label{eq:3}
	\|\lambda_S(\beta^\prime)\| \ll_{S,\eta} m(B^\prime)^{-\frac{1}{4}+\eta}\quad\hbox{for all $\eta>0$}.
	\end{equation}
	We note that this argument applies to groups that possess Cartan and Iwasawa decompositions, and it does not require that the group be simply connected.  
	Combining \eqref{eq:11}, \eqref{eq:2}, and \eqref{eq:3},
	we deduce the theorem.
\end{proof}

\subsection{Mean ergodic theorem for simply connected groups}\label{sec:mean}
In our discussion below we aim to consider approximation by elements of $\Gamma_S$ subject to arbitrary additional congruence conditions. Such a condition is determined by a compact open subset $W\subset \cI^S$, and to facilitate this discussion we will now reformulate Theorem \ref{th:norm2} in a more general and explicit form.

\begin{lem}\label{l:W}
	For a compact open subset $W\subset \cI^S$ there exists a compact open subgroup $U(W)$ of $\cI^S$
	such that $W$ is bi-invariant under $U(W)$, and $U(W)$ is a maximal subgroup of $\cU=I^S$ with this property. 
\end{lem}
\begin{proof}
	Since $W$ is open, for every $w\in W$, there exists
	a compact open subgroup $\cU_w$ of $\cI^S$ such that $\cU_w w\subset W$.
	Then by compactness, $W =\cup_{i=1}^\ell \cU_{w_i}w_i$.
	This implies that $W$ is left-invariant under the compact open subgroup
	$\cU'= \cap_{i=1}^\ell \cU_{w_i}$. A similar argument shows that
	$W$ is also right-invariant under a compact open subgroup
	$\cU''$. Therefore $U^\prime(W)=\cU'\cap \cU''$ is compact and open and leaves $W$ bi-invariant. It is clear that there exists a maximal open (and, hence, closed and compact) subgroup 
with this property, and we denote this subgroup by $U(W)$.
\end{proof}
Note that $U(W)$ depends on $S$ also, but we suppress this dependence in the notation. Let
$$
\Gamma_S(W):=\Gamma_S\cap (G_\infty \times G_S\times U(W)).
$$
Since $U(W)$ is a finite-index subgroup of $\cI^S$, it follows that 
$\Gamma_S(W)$ has finite index in $\Gamma_S$.
In particular, $\Gamma_S(W)$ is a lattice subgroup of $G_\infty\times G_S$.
We consider the homogeneous space
\begin{equation}\label{eq:def000}
Y_{S,W}:=(G_\infty\times G_S)/\Gamma_S(W)
\end{equation}
equipped with the Haar probability measure $\mu_{S,W}$.

The group $G_S$ naturally acts on the space $Y_{S,W}$ by left translations,
and we introduce averaging operators
\begin{equation}\label{eq:def001}
\pi_{S,W}(\beta):L^2(Y_{S,W})\to L^2(Y_{S,W}):\phi\mapsto 
\frac{1}{m_S(B)}\int_B \phi(g^{-1}y)\, dm_S(g).
\end{equation}
defined for measurable subsets $B$ of $G_S$ with finite positive measures.

For simplicity, we skip the index $W$ in the above notation if $W=\cI^S$.

\begin{thm} \label{th:mean_sc2}
	Assume that $\sf G$ is simply connected, and isotropic over $S\subset V_K^f$. 
	Let $W \subset \cI^S$ be a compact open subset, and let  
	$W_S$ be a compact open subgroup of $G_S$. 
	Let $\beta$ be the Haar-uniform 
	probability measure supported on a $W_S$-bi-invariant bounded subset $B$ of $G_S$ with positive measure. Then there exists $\mathfrak{k}_S({\sf G})>0$ such that
	for every $\phi\in L^2(Y_{S,W})$, and any $\eta > 0$ 
	$$
	\left\| \pi_{S,W}(\beta)\phi-\int_{Y_{S,W}}\phi\,d\mu_{S,W}\right\|_{L^2(Y_{S,W})}\ll_{\eta,W_S}
	m_{S}(B)^{-\mathfrak{k}_{ S}({\sf G})+\eta}\, \norm{\phi}_{L^2(Y_{S,W})}.
	$$
\end{thm} 

\begin{rem}\label{exponents}
{\rm 
The exponent $\mathfrak{k}_{ S}({\sf G})$ can be taken to be $\mathfrak{k}_{ S}({\sf G})=\frac{1}{4\mathfrak{n}_S({\sf G})}$ (cf.~Theorem \ref{th:norm2}). Furthermore, when the set $B$ is $U_S$-bi-invariant,
we may take 
the better exponent $\mathfrak{k}_{ S}({\sf G})=\frac{1}{\mathfrak{q}_S({\rm G})}$ (cf.~Theorem \ref{th:norm}).
}
\end{rem}

We note that the existence of the exponent $\mathfrak{k}_{ S}({\sf G})$ and its uniformity over  $W\subset \cI^S$,  is a deep property of the automorphic representation. It is a consequence of the spectral gap property of automorphic representations $\pi_v$ (explained 
above) holding uniformly over congruence subgroups. In the case of $SL_2$ for example, this property is known as the Ramanujan-Petersson-Selberg eigenvalue bounds. We refer to \cite{BS91}, \cite{BLS92}, \cite{Cl03},\cite{S05} \cite{COU}, \cite{CU} and \cite{C07} for further discussion.

\begin{proof}[Proof of Theorem \ref{th:mean_sc2}]
	To simplify notations, we set
	$$
	G:={\sf G}(\mathbb{A}_K)\quad\hbox{and}\quad\Gamma:={\sf G}(K).
	$$	
	We consider the action of the group $G_\infty\times G_S$ on the double-coset space $U(W)\backslash G/\Gamma$.
	The orbits of this action are open (and, consequently, also closed).
	Since $\sf G$ is assumed to be isotropic over $S$,
	it follows from the Strong Approximation Property (Theorem \ref{th:strong_app})
	that the projection of $\Gamma$ to $G_{V_K^f\backslash S}$ is dense.
	Therefore, the above orbits are also dense. Hence, we conclude 
	that this action is transitive, and
	$$
	Y_{S,W}\simeq U(W)\backslash G/\Gamma
	$$
	as $(G_\infty\times G_S)$-spaces.
	In particular, we also deduce equivalence
	of unitary representations
	\begin{equation}\label{eq:equiv}
	L^2(Y_{S,W}) \simeq L^2(U(W)\backslash G/\Gamma)
	\end{equation}
	of $G_S$.
	Furthermore, the space $L^2(U(W)\backslash G/\Gamma)$ can be identified with the subspace 
	$L^2(G/\Gamma)^{U(W)}$ consisting of $U(W)$-invariant
	functions in $L^2(G/\Gamma)$.
	Since $\sf G$ is assumed to be simply connected,
	there are no non-trivial automorphic characters (see, for instance, \cite[Lem.~4.1]{GGN1}). Hence, it follows from Theorem \ref{th:norm2}
	(or from Theorem \ref{th:norm} if $B$ is $U_S$-bi-invariant)
	that there exists $\mathfrak{k}_S({\sf G})>0$ such that
	for every $\psi\in L^2(G/\Gamma)$, and any $\eta > 0$ 
	$$
	\left\| \pi^{\rm aut}_{S}(\beta)\psi-\mathcal{P}(\psi)\right\|_{L^2(G/\Gamma)}\ll_{W(S),\eta}
	m_{S}(B)^{-\mathfrak{k}_{ S}({\sf G})+\eta}\, \norm{\psi}_{L^2(G/\Gamma)},
	$$
	where $\mathcal{P}$ denotes the orthogonal projection on the space of constant functions. Hence, the statement of the theorem follows from \eqref{eq:equiv}.
\end{proof}

We shall also prove and utilize a version of Theorem \ref{th:mean_sc2}
for general $K$-simple groups, but we will postpone this until \S \ref{sec:mean2} below.

\subsection{Riemannian local volume and distance estimates}
We now note two local properties of distance and volume in any almost connected  semisimple Lie group with finite center, denoted by $G_\infty$, which will be used in our arguments below. 
We fix a right invariant Riemannian metric $\rho$ on $G_\infty$ and consider the corresponding balls
$$
B(g,r):=\{x\in G_\infty:\, \rho(x,g)\le r \}.
$$
We first establish the following local estimates for the volume
$m_\infty(B(g,r))$. By right-invariant, it clearly suffices to consider only the case where the center is the identity element $e$. Recall that we denote $d=\dim_\RR G_\infty$.

\begin{lem}\label{l:ball}
\begin{enumerate}
\item[(a)] There exist $c_1,c_2>0$ and $r_0>0$ such that
$$
c_1\, r^d\le m_\infty(B(e,r)) \le c_2\, r^d\quad\hbox{for all $r\in (0,r_0)$.}
$$
\item[(b)] For $0< c_0 < 1$, there exist $c,r'_0>0$ such that for all $\epsilon\in (0,c_0\,r)$ and $r\in (0,r'_0)$,
$$
m_\infty(B(e,r+\epsilon))-m_\infty(B(e,r))\le c\,\frac{\epsilon}{r}\, m_\infty(B(e,r)).
$$	
\end{enumerate}
\end{lem}

\begin{proof}
According to the volume formula for Riemannian balls \cite[p.~66]{Sak},
for sufficiently small $r$,
$$
m_\infty(B(e,r))=\int_0^r \omega(s)\, ds,
$$
where $\omega$ is a continuous function satisfying
$$
c_1'\, s^{d-1}\le \omega(s)\le c_2'\, s^{d-1}
$$
for some $c_1',c_2'>0$. This implies the first estimate. Also,
$$
m_\infty(B(e,r+\epsilon))-m_\infty(B(e,r))=\int_r^{r+\epsilon}\omega(s)\, ds
\ll (r+\epsilon)^{d-1}\, \epsilon, 
$$
which gives the second bound.
\end{proof}	
Let us now define two metrics $\rho$ and $\rho^\prime$ on $G_\infty$ to be locally equivalent if for every compact neighborhood $Q\subset G_\infty$ of $I$, there exists a constant $C_Q$ such that for any two points $x\neq y\in G$ satisfying $xy^{-1}\in Q$, we have  
$C_Q^{-1}\le \frac{\rho(x,y)}{\rho^\prime(x,y)}\le C_Q$. 

Assuming that $G_\infty\subset GL_d(\RR)$ (namely fixing a faithful linear represesentation of $G_\infty$, but suppressing it from the notation), fix a (vector space) norm on $\text{Mat}_d(\RR)$, and consider the distance on $G_\infty$ given by $\rho^\prime(x,y)=\norm{x-y}$, $x,y \in G_\infty$. 

\begin{lem}\label{l:norm-equiv}
Any right (or left) invariant Riemannian distance $\rho$ on $G_\infty$ is locally equivalent to the distance $\rho^\prime$ defined by any (vector space) norm in any faithful linear representation (as defined above).  
\end{lem}

\begin{proof}
The right-invariant Riemannian metric on $G_\infty$ defining the right-invariant distance $\rho$ is determined by the choice of a positive-definite inner product on the Lie algebra $\mathfrak{g}_\infty$ of $G_\infty$. The inner product  determines a Euclidean norm $\abs{X}$ on the Lie algebra. We have 
  $\mathfrak{g}_\infty\subset \text{Mat}_d(\RR)$, and the exponential map denoted $X\to e^X$ takes $\mathfrak{g}_\infty$ into $G$ and is a diffeomorphism on a (Euclidean) ball $D_{\vare_0}(0)$ centered at $0\in\mathfrak{g}_\infty$. For $X\in D_{\vare_0}(0)$ we denote $x=e^X\in G$, and then $\rho(I, x)=\rho(I, e^X)=\abs{X}+O(\abs{X}^2)$. In $\text{Mat}_d(\RR)$ we have $\norm{I-x}=\norm{I-e^X}=\norm{X}+O(\norm{X}^2)$ when $\norm{X} < 1/2$ (say). 
Hence there exists $C_{Q_0}$ satisfying that $C^{-1}_{Q_0}\le \rho(I,x)/\norm{I-x}\le C_{Q_0}$ for $x\neq I$ in the compact neighborhood of $I$ given by the set $ Q_0=\set{e^X\,;\,X\in D_{\vare_1}(0)}\subset G_\infty$, for suitable $\vare_1 >0$. 

Since the operator norms of $y$ and $y^{-1}$ are bounded above and below when $y\in Q_0$, writing  $\norm{x-y}=\norm{(I-xy^{-1})y}\le \norm{y}_{op}\norm{I-xy^{-1}}$, and 
$\norm{x-y}\ge \norm{I-xy^{-1}}/\norm{y^{-1}}_{op}$, we 
conclude that $C^{-1}_{Q_1}\le \frac{\rho(y,x)}{\norm{y-x}}\le C_{Q_1}$ for any $x\neq y$ in compact neighborhood $Q_1\subset Q_0$ of $I$. It then follows by compactness and continuity that the same holds (with a different constant $C_{Q_R}$) for $x\neq y$ in any compact neighborhood $Q_{R}$ of $I$.  

Now fixing any compact set $Q$, given any $x\in G_\infty$ the set $Qx$ is contained in compact neighborhood $Q_R$ of $I$ for some $R$ (depending on $x$), and so for all $y\in Qx, y\neq x$ we have 
$C^{-1}_{Q}\le \frac{\rho(y,x)}{\norm{y-x}}\le C_{Q}$, as stated. 
\end{proof}

\section{Asymptotic formula for the counting function of Diophantine approximants} 
\label{sec:schmidt}

%

Let $\sf G$ be a simply connected $K$-simple linear algebraic group defined over a number field $K$. We will freely use the notation introduced in Sections \ref{sec:intro}--\ref{sec:not}. 
Our goal in the present section is to establish an analogue of W. Schmidt's theorem. We will analyze the number of Diophantine approximants
$r\in \Gamma_S$ to points $x\in G_\infty$, collect them into the  
counting function :
\begin{equation}\label{NT}
N_T(x):=\big|\{r\in \Gamma_S:\, \rho(x,r)\le \h (r)^{-b},\; 1\le \h(r)\le T\}\big|.
\end{equation}
and establish its asymptotics. 
We refer to the exponent $b$ as the scale of approximation. 

\subsection{The parameters of effective Diophantine approximation}
Fix $S\subset V_K^f$.
We recall that the group $\Gamma_S$ embeds diagonally in $G_\infty\times G_S$
as a lattice subgroup, and the Haar measures $m_\infty$ and $m_S$
on the factors are normalized so that 
$\Gamma_S$ has covolume one with respect to $m_\infty\times m_S$.
We expect that $N_T(x)$ is approximated by the volume sum
\begin{equation}\label{eq:vvv}
V_T:={\sum}_{1\le h\le T} m_\infty(\cB_h)m_S(\Sigma_S(h)),
\end{equation}
where 
$$
\cB_h:=\left\{x\in G_\infty:\, \rho(x,e)\le h^{-b}\right\}\quad\hbox{and}\quad \Sigma_S(h):=\{g\in G_S:\, \h(g)=h\}.
$$ 
We recall that $\cI_S=\prod_{v\in S} {\sf G}(O_v)$.
Clearly, the set $\Sigma_S(h)$ is compact and $\cI_S$-bi-invariant,
but it might be empty for some choices of $h\in \N$,
so that we introduce
$$
\mathcal{L}_S:=\{h\in\N:\, \Sigma_S(h)\ne \emptyset\}.
$$
The estimates in this section will depend on the following three parameters:
\begin{center}
	\begin{tabular}{ll}
$\mathfrak{a}$ &-- the volume growth rate of the sets $\Sigma_S(h)$ (see \eqref{sphere-vol}),\\
 $\mathfrak{k}$ &-- the error term in the Ergodic Theorem
(see \eqref{sphere-norm}),\\
$\mathfrak{d}$ &-- the volume decay rate  of balls $\mathcal{B}_h$ (see \eqref{vol-II}).
	\end{tabular}
\end{center}

We record some basic properties of the sets $\Sigma_S(h)$ and $\mathcal{L}_S$:

\begin{lem}\label{l:sigma}
Let $S$ be a finite subset of $V_K^f$.
	
	\begin{enumerate}
		\item[(a)] For every $\eta>0$,  
		$$
		{\sum}_{T\in \mathcal{L}_S} T^{-\eta}<\infty\quad\hbox{and}\quad \big|\mathcal{L}_S\cap [1,T]\big|=O_\eta(T^\eta).
		$$
		\item[(b)] The set $\log (\mathcal{L}_S)$ has bounded gaps.
		\item[(c)] 
		Suppose that ${\sf G}$
		is isotropic over $K_v$ for all $v\in S$. 
		Then there exists an exponent $\mathfrak{a}=\mathfrak{a}(S)>0$ such that
		\begin{equation}\label{sphere-vol}\tag{V1}
		m_{S}(\Sigma_S(h))\gg_S h^{\mathfrak{a}}\quad\quad \hbox{for all $h\in \mathcal{L}_S$.}
		\end{equation}
	\end{enumerate}
\end{lem}

We note that in order to obtain the lower bound in terms of $h$ in part (c), it is essential to assume that the group $\sf G$ is isotropic for all $v\in S$.
This is the only place where we use this condition. When this condition does not hold, it is still possible to formulate a weaker lower bound. 

\begin{proof}[Proof of Lemma \ref{l:sigma}]
	We observe that for $g\in G_S$, all values of the height $\h(g)$ are of the form $\prod_{v\in S} q_v^{n_v}$, where $q_v$ denotes the norm of uniformizing parameter of $K_v$.	Since $S$ is assumed to be finite, this implies the claim (a). 
	
	To prove (b), it will be convenient to consider $\sf G$ as a subgroup of $\hbox{SL}_n$. We note that when $S_1\subset S_2$, we have $\mathcal{L}_{S_1}\subset \mathcal{L}_{S_2}$. 
	Hence it is sufficient to prove claim (b) when $S$ consists of a single place $v$. 
	Since $\sf G$ is isotropic over this place,
	it contains a non-trivial $K_v$-split torus $\sf A$. 
	There exists $g\in \hbox{SL}_n(K_v)$ such that $g{\sf A}g^{-1}$ is diagonal.
	Let us fix $a\in {\sf A}(K_v)$ such that $\|gag^{-1}\|_v>1$. 
	Since $\|ga^ng^{-1}\|_v=\|gag^{-1}\|^n_v$, it is clear that the set 
	$\left\{\log \|ga^ng^{-1}\|_v:\, n\in\mathbb{N}\right\}$ has bounded gaps. Finally, we note that for some $c>1$, $c^{-1}\|z\|_v\le \|gzg^{-1}\|_v\le c \|z\|_v$ for all $z\in \hbox{Mat}_n(K_v)$.
	This implies that the set $\{\log \|a^n\|_v:\, n\in\mathbb{N}\}$ has bounded gaps as well.
	
	Claim (c) was established in \cite[Lemma~4.3]{gk} in the case of the field of rationals, and this argument generalises to general number fields. 
\end{proof}

For $h\in \mathcal{L}_S$, we denote by $\sigma_h$ the Haar-uniform probability measure supported on the subset $\Sigma_S(h)^{-1}$ of $G_S$ 
and consider the corresponding averaging operator
$$
\pi_S(\sigma_h):L^2(Y_S)\to L^2(Y_S)
$$
on the space $Y_S:=(G_\infty\times G_S)/\Gamma_S$ defined in \eqref{eq:def001}.
We note that $\Sigma_S(h)^{-1}\ne \Sigma_S(h)$ in general, but both sets have the same Haar measure.
By Theorem \ref{th:mean_sc2},
there exists $\mathfrak{k}_S({\sf G})=\mathfrak{k}_S\in (0,1/2)$ such that for all $\phi\in L^2(Y_S)$, and any $\eta > 0$ 
\begin{equation}\label{sphere-norm}\tag{SP}
\norm{\pi_S(\sigma_h)\phi -\int_{Y_S} \phi\, d\mu_S}_{L^2(Y_S)}\ll_{,\eta} \,  m_{S}(\Sigma_S(h))^{-\mathfrak{k}_S+\eta}\quad\hbox{for all $h\in \mathcal{L}_S$.} 
\end{equation}

We also recall that by Lemma \ref{l:ball}(a),
\begin{equation}\label{vol-II}\tag{V2}
h^{-b\mathfrak{d}}\ll m_{\infty}(\mathcal{B}_h)\ll h^{-b\mathfrak{d}}\quad\hbox{for all $h\ge h_0$}
\end{equation}
where $\mathfrak{d}:=\dim_\R(G_\infty)$.

We set 
$$
b_0:=2\mathfrak{a}\mathfrak{k}/\mathfrak{d}
$$
and 
for a positive scale  $b<b_0$, define 
$$
\theta_0(b):=\frac{(1-\mathfrak{k})\mathfrak{a}-b\mathfrak{d}/2}{\mathfrak{a}-b\mathfrak{d}}=\frac{1}{2}+\frac{(1/2-\mathfrak{k})\mathfrak{a}}{\mathfrak{a}-b\mathfrak{d}}.
$$
Note that the condition on $b<b_0$ insures that $\theta_0(b)\in (0,1)$.
With these notations, we prove:

\begin{thm}\label{ae-disc}
Let $\sf G$ be a connected simply-connected $K$-simple algebraic group defined over a number field $K$, and $S$ is a finite set of finite places
such that $\sf G$ is isotropic for all $v\in S$. 
Then for every $b\in (0,b_0)$ and $\theta\in (\theta_0(b), 1)$, 
$$
\big\|N_T-V_T\big\|_{L^2(Q)}\ll_{S, Q,\theta} V_T^\theta,
$$
where $Q$ is an arbitrary bounded measurable subset of $G_\infty$.
\end{thm}

The proof of Theorem \ref{ae-disc} will be based on the estimates (V1), (V2), (SP).

\vspace{0.2cm}

Using a Borel--Cantelli argument we will also derive a pointwise bound for $N_T$:

\begin{cor}\label{cor_NT}
With the notation of Theorem \ref{ae-disc}, for every $\theta'>\theta$
$$
N_T(x)=V_T+O_{S, x,\theta'}\big(V_T^{\theta'}\big)\quad\hbox{for a.e. $x\in G_\infty$.}
$$ 
\end{cor}

\subsection{ Towards a best-possible estimate.}\label{sec:best}
Let us note that it is often the case that the lower bound $m_S(\Sigma_S(h))\ll_S h^{\mathfrak{a}}$ holds, as well as the upper bound (\ref{sphere-vol}). Then the sum \eqref{eq:vvv} is bounded by
$$	 
V_T\ll {\sum}_{1\le h\le T,h\in\mathcal{L}_S}  h^{-b\mathfrak{d}} \cdot h^{\mathfrak{a}}\,,	
$$
and it follows from Lemma \ref{l:sigma}(a) that $V_T$
is uniformly bounded when $b>\mathfrak{a}/\mathfrak{d}$. Therefore, 
the estimate in Theorem \ref{ae-disc} is only interesting 
in the range $b\le\mathfrak{a}/\mathfrak{d}$.
We highlight that Corollary \ref{cor_NT} gives the following estimate :

\begin{cor}\label{rem:temp}
Suppose that additionally in Corollary \ref{cor_NT}	
the automorphic representation of $G_S$ is tempered 
and $S$ consists of unramified places. Then for every 
$b<\mathfrak{a}/\mathfrak{d}$,
$$
N_T(x)=V_T+O_{S, x,\eta}\left(V_T^{\frac{1}{2}+\eta}\right)\quad\hbox{for all $\eta>0$ and a.e. $x\in G_\infty$.}
$$ 
\end{cor}

Indeed, under the temperedness condition the estimate \eqref{sphere-norm}
holds with $\mathfrak{k}_S=1/2$ (see Theorem \ref{th:norm}), so that in this case 
Theorem \ref{ae-disc} covers all the relevant range of parameters $b<\mathfrak{a}/\mathfrak{d}$ and gives the exponent
$\theta=\frac{1}{2}+\eta$ for every $\eta>0$. Hence,
Corollary \ref{rem:temp} is a direct consequence of Corollary \ref{cor_NT}.

Noting the fact that error term in the foregoing estimate is bounded (in essence) by the square root of the main term, it is natural to expect that the exponent it establishes is in fact best-possible. A full proof of this fact requires establishing the expected lower bound for $N_T$. We will establish lower bounds for discrepancy estimates in a separate paper, but the question of optimality of Corollary \ref{rem:temp} remains open. 
\subsection{Proof of Theorem \ref{ae-disc}} 
Consider the function
$$
D_T(x):=N_T(x)-V_T,\quad \hbox{with $x\in G_\infty$}\,,
$$
and recall that the metric $\rho$ is right-invariant, and so $B(x,h^{-b})= \cB_{h}x= \cB_{h}^{-1}x$. Therefore, 
\begin{align*}
D_T(x):=&\sum_{1\le h\le T} \Big(\big|\{\gamma \in\Gamma_S\cap B(x,h^{-b}):\, \h(\gamma)=h\}\big|-m_{\infty}(\cB_{h}) m_{S}\big(\Sigma_S(h)\big)\Big)\\
=& \sum_{1\le h\le T} \Big(\big|\Gamma_S\cap (\cB_{h}x\times \Sigma_S(h))\big|-m_{\infty}(\cB_{h}) m_{S}\big(\Sigma_S(h)\big)\Big).
\end{align*}
The crucial ingredient of our proof is that $D_T(x)$ can be 
represented in terms of the averaging operators $\pi_S(\sigma_h)$.
Let 
$$
\chi_h(g_\infty,g_S):=\chi_{\cB_h}(g_\infty)\chi_{\cI_S}(g_S)
$$
denote the characteristic function of the subset 
$\cB_h \times \cI_S$ of $G_\infty\times G_S$.
The sum 
$$
\phi_h (g):=\sum_{\gamma\in \Gamma_S} \chi_h(g\gamma)=\sum_{\gamma\in \Gamma_S} \chi_h(g\gamma^{-1})
$$
defines a measurable function with compact support
on the homogeneous space $Y_S:=(G_\infty\times G_S)/\Gamma_S$. 
We observe that for $x\in G_\infty$ and $u\in \cI_S$,
\begin{align*}
\int_{a\in \Sigma_S(h)} \phi_h(a (x,u))\,dm_{S}(a)
&= \sum_{\gamma\in \Gamma_S} \int_{\Sigma_S(h)} \chi_h(x\gamma^{-1}, au\gamma^{-1})\, d m_{S}(a)\\
&= \sum_{\gamma\in \Gamma_S} \int_{\Sigma_S(h)} \chi_{\cB_h}( x\gamma^{-1}) \chi_{\cI_S}(au\gamma^{-1})\, dm_{S}(a)\\
&= \sum_{\gamma\in \Gamma_S\cap B(x,h^{-b})} m_{S}\big(\cI_S\gamma u^{-1}\cap \Sigma_S(h)\big).
\end{align*}
Note that in the above computation, since $\Gamma_S$ is a discrete subgroup of $G_\infty\times G_S$, the non-zero summands in 
the sum above  constitute a finite subset of $\gamma$'s. 
In order to evaluate the last expression, we use that the set 
$\Sigma_S(h)$
is $\cI_S$-bi-invariant. Therefore,
if $H_f(\gamma)=h$, then $ \cI_S \gamma u^{-1}\subset \Sigma_S(h)$ and $m_{S}( \cI_S\gamma u^{-1}\cap \Sigma_S(h))$=1.
On the other hand, if $H_f(\gamma)\neq h$, we have $\cI_S\gamma u^{-1}\cap \Sigma_S(h)=\emptyset$.
Hence, for every $u\in \cI_S$,
$$
\big|\{\gamma \in\Gamma_S\cap B(x,h^{-b}):\, \h(\gamma)=h\}\big|=\int_{a\in \Sigma_S(h)} \phi_h(a (x,u))\,dm_{S}(a)\,.
$$
Furthermore, 
\begin{align*}
\int_{Y_S} \phi_h \, d\mu_S &=\int_{G_\infty\times G_S} \chi_{\cB_h}(g_\infty)\, \chi_{\cI_S}(g_S)\, dm_\infty(g_\infty) d m_S(g_S)\\
&= m_{\infty}(\cB_h)m_S(\cI_S)= m_{\infty}(\cB_h).
\end{align*}
Hence, we deduce that for every $u\in \cI_S$,
$$
D_T(x)=\sum_{1\le h\le T} \left(\int_{a\in \Sigma_S(h)} \phi_{h}(a(x,u))\,d m_{S}(a)
- m_{S}\big(\Sigma_S(h)\big)\int_{Y_S} \phi_h \, d\mu_S\right).
$$
Let ${Q}$ be a bounded measurable subset in $G_\infty$. Since the first integral in the previous line was just shown to be independent of  $u\in \cI_S$, we obtain
\begin{align*}
&\int_{{Q}} \left|\int_{a\in \Sigma_S(h)} \phi_h(a (x,e))\,d m_{S}(a)
- m_S(\Sigma_S(h))\int_{Y_S} \phi_h \, d\mu_S\right|^2\, d m_{\infty}(x)\\
=&\int_{{Q}\times \cI_S} \left|\int_{a\in \Sigma_S(h)} \phi_h(a (x,u)\,d m_{S}(a)
- m_{S}(\Sigma_S(h))\int_{Y_S} \phi_h \, d\mu_S\right|^2\, d m_{\infty}(x)d m_{S}(u).
\end{align*}
The set $Q\times \cI_S$ projects onto a measurable subset $(Q\times \cI_S)\Gamma_S$ of $Y_S$.
Since $Q\times \cI_S$ is bounded, there exists $N_Q$ such that
every point in the image has preimage of cardinality at most $N_Q$.
This implies that the last integral is bounded by
\begin{align*}
&N_Q \int_{({Q}\times \cI_S)\Gamma_S} \left|\int_{a\in \Sigma_S(h)} \phi_h(a y)\,d m_{S}(a)
- m_{S}(\Sigma_S(h))\int_{Y_S} \phi_h \, d\mu_S\right|^2\, d\mu_S(y)\\
\le&\, N_{{Q}}\int_{Y_S} \left|\int_{\Sigma_S(h)} \phi_h(a y)\,d m_{S}(a)
- m_{S}(\Sigma_S(h))\int_{Y_S} \phi_h \, d\mu_S\right|^2\, d\mu_S(y)\\
=&\, N_{{Q}}\, m_{S}(\Sigma_S(h))^2\left\|\pi_S(\sigma_h)\phi_h -\int_{Y_S} \phi_h \, d\mu_S \right\|_{L^2(Y_S)}^2.
\end{align*}
Hence, using the bound (\ref{sphere-norm}), we deduce that for any $\eta > 0$ 
\begin{align*}
\|D_T\|_{L^2(Q)} &\le \sum_{1\le h\le T} \left\|\int_{a\in \Sigma_S(h)} \phi_{h}(a (\cdot,e))\,dm_{S}(a)
-m_{S}(\Sigma_S(h))\int_{Y_S} \phi_h \, d\mu_S\right\|_{L^2(Q)}\\
& \ll_{S,Q,\eta}   \sum_{1\le h\le T} m_{S}(\Sigma_S(h))^{1-\mathfrak{k}_S+\eta}\|\phi_h\|_{L^2(Y_S)}.
\end{align*}
Furthermore, the $L^2$-norm of $\phi_h$ can be estimated as follows
\begin{align*}
\|\phi_h\|^2_{L^2(Y_S)} &=\int_{(G_\infty\times G_S)/\Gamma_S} \left(\sum_{\gamma\in \Gamma_S} \chi_h(g\gamma)\right)^2\,d(m_\infty\times m_S)(g)\\
&=\int_{(G_\infty\times G_S)/\Gamma_S} \sum_{\gamma_1,\gamma_2\in \Gamma_S} \chi_h(g\gamma_1)\chi_h(g\gamma_2)\,d(m_\infty\times m_S)(g)\\
&=\int_{(G_\infty\times G_S)/\Gamma_S} \sum_{\gamma,\delta\in \Gamma_S} \chi_h(g\delta)\chi_h(g\delta\gamma)\,d(m_\infty\times m_S)(g)\\
&= \int_{G_\infty\times G_S} \sum_{\gamma\in \Gamma_S} \chi_h(g)\chi_h(g\gamma)\,d(m_\infty\times m_S)(g)\\
& =\sum_{\gamma\in \Gamma_S} (m_\infty\times m_S)\left((\cB_h\times \cI_S)\cap (\cB_h \times \cI_S)\gamma^{-1}\right)\le 
N_h\, m_{\infty}(\cB_h)\,,
\end{align*}
where $N_h$ denotes the number of $\gamma\in \Gamma_S$ such that 
$(\cB_h\times \cI_S)\cap (\cB_h \times \cI_S)\gamma^{-1}\ne \emptyset$.
Since the family of sets $\cB_h\times \cI_S$ is uniformly bounded,
this number is uniformly bounded, and we conclude using \eqref{vol-II} that for any $\eta > 0$ 
$$
\|D_T\|_{L^2(Q)}\ll_{S,Q,\eta} \sum_{1\le h\le T} m_{S}(\Sigma_S(h))^{1-\mathfrak{k}_S+\eta}  m_{\infty}(\cB_h)^{1/2}\ll
\sum_{1\le h\le T} m_{S}(\Sigma_S(h))^{1-\mathfrak{k}_S+\eta} h^{-b\mathfrak{d}/2}.
$$
To complete the proof of Theorem \ref{ae-disc} we employ the following computation, where we use parameters $A>0$ and $B>1$ that will be specified later. Using \eqref{sphere-vol}, and writing $\mathfrak{k}_S-\eta=\mathfrak{k}$ for brevity, we obtain
\begin{align*}
\sum_{1\le h\le T} m_{S}(\Sigma_S(h))^{1-\mathfrak{k}} h^{-b\mathfrak{d}/2}
\ll_S \sum_{1\le h\le T} m_{S}(\Sigma_S(h))^{1-\mathfrak{k}+A} h^{-(\mathfrak{a}A+b\mathfrak{d}/2)},
\end{align*}
and applying H\"older's inequality, we conclude that this sum is bounded by
$$
\left(\sum_{1\le h\le T} m_{S}(\Sigma_S(h))^{(1-\mathfrak{k}+A)B} h^{-(\mathfrak{a}A+b\mathfrak{d}/2)B}\right)^{1/B}\big|\mathcal{L}_S\cap [1,T]\big|^{1-1/B}.
$$
We choose $A$ and $B$ so that
$$
(1-\mathfrak{k}+A)B=1\quad\hbox{and}\quad (\mathfrak{a}A+b\mathfrak{d}/2)B=b\mathfrak{d},
$$
namely,
$$
A=\frac{(1-\mathfrak{k})b\mathfrak{d}-b\mathfrak{d}/2}{\mathfrak{a}-b\mathfrak{d}}\quad\hbox{and}\quad
B=\frac{\mathfrak{a}-b\mathfrak{d}}{(1-\mathfrak{k})\mathfrak{a}-b\mathfrak{d}/2}.
$$
Taking into account that $\mathfrak{k}\in (0,1/2)$ and $b<2\mathfrak{k}\mathfrak{a}/\mathfrak{d}$,
a direct computation verifies that $A>0$ and $B>1$.
Hence, using Lemma \ref{l:sigma}(a) and estimate \eqref{vol-II},
we conclude that  for every $\eta>0$,
$$
\|D_T\|_{L^2(Q)}\ll_{S,Q,\eta} 
\left(\sum_{1\le h\le T} m_{S}(\Sigma_S(h)) h^{-b\mathfrak{d}}\right)^{1/B} T^{\eta(1-1/B)}\ll_\eta V_T^{1/B} T^{\eta(1-1/B)}.
$$
Finally, we note that it follows from \eqref{sphere-vol}--\eqref{vol-II} and Lemma \ref{l:sigma}(b) that 
\begin{equation}\label{eq:lower_b}
V_T\gg_S \sum_{h\in\mathcal{L}_S\cap [1,T]} h^{-b\mathfrak{d}}\cdot h^{\mathfrak{a}}\gg_S T^{\mathfrak{a}-b\mathfrak{d}}.
\end{equation}
Since $\mathfrak{a}-b\mathfrak{d}>0$, we deduce from the previous estimate that for every $\eta'>0$,
$$
\|D_T\|_{L^2(Q)}\ll_{Q,\eta'} V_T^{\theta_0+\eta'},
$$
where 
$$
\theta_0:=1/B=\frac{(1-\mathfrak{k})\mathfrak{a}-b\mathfrak{d}/2}{\mathfrak{a}-b\mathfrak{d}}.
$$ Since $B>1$, we have $\theta_0\in (0,1)$.
This completes the proof of Theorem  \ref{ae-disc}. \qed

\begin{proof}[Proof of Corollary \ref{cor_NT}]
We note that the function 
$D_T(x)=N_T(x)-V_T$
is piecewise constant in $T$ and is determined by
its values for $T\in \mathcal{L}_S$. Hence, it is sufficient to analyse
the values $D_T(x)$ for $T\in \mathcal{L}_S$.
Given any bounded measurable subset $Q$ of $G_\infty$,
we have established in Theorem \ref{ae-disc} the bound
$$
\int_{{Q}} |D_T(x)|^2\, dm_{\infty}(x)\ll_{S,{Q},\theta} V_T^{2\theta}.
$$
Therefore, for every $\eta>0$, the sets
$$
\Omega_T:=\{x\in {Q}:\, |D_T(x)|\ge  V_T^{\theta}\,T^\eta \}
$$
satisfy
$$
m_{\infty}(\Omega_T)\ll_{S,{Q},\theta} T^{-2\eta}.
$$
Hence, it follows from Lemma \ref{l:sigma}(a) that
$$
\sum_{T\in\mathcal{L}_S} m_{\infty}(\Omega_T)<\infty,
$$
and by the Borel--Cantelli Lemma, the $\limsup$ of the sets $\Omega_T$
with $T\in\mathcal{L}_S$ has measure zero. This means that for almost every $x\in {Q}$ and $T\in \mathcal{L}_S$, 
$$
|D_T(x)|\le V_T^{\theta}\, T^{\eta}\quad\hbox{ when $T\ge T_0(x,\eta)$}.
$$
Then using the estimate \eqref{eq:lower_b}, we conclude that
for almost all $x\in Q$,
$$
|D_T(x)|\ll_S V_T^{\theta+\eta/(\mathfrak{a}-b\mathfrak{d})}
\quad\hbox{ when $T\ge T_0(x,\eta)$}.
$$
Since $G_\infty$ can be exhausted by a countable union of bounded measurable sets $Q$, it follows that for every $\theta'>\theta$
and almost all $x\in G_\infty$,
$$
|D_T(x)|\ll_S V_T^{\theta'}\quad\hbox{ when $T\ge T_0(x,\theta')$}.
$$
Thus, in particular,
$$
|D_T(x)|\ll_{S,x,\theta'} V_T^{\theta'}\quad \hbox{for all $T$,} 
$$
which implies the corollary.
\end{proof}

\section{Discrepancy bounds for simply connected groups}\label{sec:dis_simp}

Let $\sf G$ be a connected $K$-simple algebraic group defined over a number field $K$. In the present section, we assume that $\sf G$ is simply connected, and establish three effective discrepancy estimates for the distribution of rational points, namely mean-square, almost sure, and pointwise everywhere estimates, for general sets $E\subset G_\infty$.

\subsection{Discrepancy of rational points}
We fix a subset $S$ of finite places of $K$ such that $\sf G$ is isotropic over $S$ and consider the $S$-arithmetic group $\Gamma_S:={\sf G}(O_S)$
which is exhausted by the increasing family of subsets
$$
R_S(h):=\{\gamma\in \Gamma_S:\, \h(\gamma)\le h \}.
$$
We recall that $\cI_S:=\prod_{v\in S} {\sf G}(O_v)$,
and $\cI^S:=\prod_{v\in V_K^f\backslash S} {\sf G}(O_v)$, and
the Haar measure $m_S$ on $G_S$ is normalized so that $m_S(\cI_S)=1$.
We also denote by $m^S$ the Haar measure on $\cI^S$ such that $m^S(\cI^S)=1$.
Under the diagonal embedding 
$$
\Gamma_S \hookrightarrow G_\infty\times G_S,
$$
$\Gamma_S$ is a lattice subgroup, and so the sets $R_S(h)$, while infinite, are locally finite namely deposit in every bounded subset of $G_\infty$ only finitely many elements.
The Haar measure $m_\infty$ in $G_\infty$ is normalized 
so that $\Gamma_S$ has covolume one with respect to $m_\infty\times m_S$.

By the Strong Approximation Property (Theorem \ref{th:strong_app}),
the diagonal embedding 
$$
\Gamma_S \hookrightarrow G_\infty\times \cI^S
$$
is dense. Our goal is to analyse the discrepancy of distribution of this
dense set. One can show (in fact, it follows from our results here)
that the number of points from $R_S(h)$
contained in a bounded subset $\Omega$  of $G_\infty\times \cI^S$ grows as 
$$
v_S(h):=m_S\big(B_S(h)\big),
$$
where $B_S(h):=\{g\in G_S:\, \h(g)\le h \}$.
We define the discrepancy function as 
\begin{equation}\label{disc-def-2}
\mathcal{D}(R_S(h),\Omega):= \left|\frac{|R_S(h)\cap \Omega|}{v_S(h)}-(m_\infty\times m^S)(\Omega)\right|.
\end{equation}
Our goal is to produce an explicit estimate for this quantity 
for a natural collection of subsets of $G_\infty\times \cI^S$.
Let $E$ be a subset of $G_\infty$ and $W$ a compact open subset of $\cI^S$.
For $x\in G_\infty$, we set 
$$
E(x):=Ex.
$$
We will focus on analyzing the discrepancy for $\Omega=E(x)\times W$, where we allow an arbitrary congruence condition $W\subset \cI^S$.

A crucial ingredient of our analysis is the estimate on averaging operators established in \S \ref{sec:mean}, which  we now recall. We consider the spaces 
$$
Y_{S,W}:=(G_\infty\times G_S)/\Gamma_S(W)
$$
equipped with the invariant probability measures $\mu_{S,W}$.
Let $\beta_h$ be the uniform probability measure supported on 
the set $B_S(h)^{-1}$ and 
$$
\pi_{S,W}(\beta_h):L^2(Y_{S,W})\to L^2(Y_{S,W})
$$
the corresponding averaging operator defined in \eqref{eq:def001}.
We note that $B_S(h)^{-1}\ne B_S(h)$ in general.
According to Theorem~\ref{th:mean_sc2}, there exists an exponent $\mathfrak{k}_S({\sf G})=\mathfrak{k}_S\in (0,1/2]$, which is uniform in $W$, such that
for every $\phi \in L^2(Y_{S,W})$, and every $\eta > 0$ 
\begin{equation}\label{eq:nnorm} \tag{SP}
\left\| \pi_{S,W}(\beta_h)\phi-\int_{Y_{S,W}}\phi\,d\mu_{S,W}\right\|_{L^2(Y_{S,W})}\ll_{W,S,\eta}
m_{S}\big(B_S(h)\big)^{-\mathfrak{k}_S+\eta}\, \norm{\phi}_{L^2(Y_{S,W})}.
\end{equation}

The estimates in this section will depend on the parameters:

\begin{center}
\begin{tabular}{ll}
$\mathfrak{k}_S$ & -- the speed of convergence in the effective mean ergodic theorem in \eqref{eq:nnorm},\\
$\mathfrak{d}$ & -- $\dim_\RR G_\infty$, namely the decay rate of volume of balls of small radius in \eqref{vol-II}.
\end{tabular}
\end{center}

\subsection{Mean-square discrepancy estimates}
We will begin our discussion by establishing an $L^2$-bound for the discrepancy function $\mathcal{D}(R_S(h),E(x)\times W)$. Remarkably, this bound holds in great generality
for measurable subsets $E$, and as noted already it is uniform in $W$.

For a subset $E$ of $G_\infty$, we define  
$$
\mathcal{N}(E):=|\{\gamma\in {\sf G}(O):\, m_\infty(E\cap \gamma E)>0 \}|<\infty.
$$
Clearly, $\mathcal{N}(E)$ is uniformly bounded when $G_\infty$ is compact, and finite if $E$ is bounded, 
but this is not the case in general.

 \begin{thm}[mean square discrepancy bound] \label{th:l2_bound}
 Let $E$ be any measurable subset of $G_\infty$ of positive finite measure 
 satisfying $\mathcal{N}(E)<\infty$
 and $W$ a compact open subset of $\cI^S$. Then for any $\eta > 0$ 
	$$
	\Big\|\mathcal{D}(R_S(h), E(\cdot)\times W)\Big\|_{L^2(Q)}\ll_{S,Q, \eta} \mathcal{N}(E)^{1/2}
	m_\infty(E)^{1/2}m^S(W)^{1/2} v_S(h)^{-\mathfrak{k}_S({\sf G})+\eta}
	$$ 
for every bounded measurable subset $Q$ of $G_\infty$.
\end{thm}

Note that in foregoing formula, the discrepancy decreases as the size of $W$ decreases, and this fact reflects our choice of normalization in (\ref{disc-def-2}).  If we choose to normalize the discrepancy of rational points in $E$ by dividing by the measure of $W$, the error estimate on the r.h.s. would depend on $W$ via  $m^S(W)^{-1/2}$. Then 
the discrepancy will increase as the congruence conditions imposed become more stringent, as expected. 

\begin{rem}
{\rm 
We say a subset ${Q}$ of $G_\infty$ is 
\emph{$\Gamma_S(W)$-injective} if the 
projection map
$$
Q\times \cI_S \to \left(G_\infty\times G_S\right)/\Gamma_S(W)
$$
is injective.
If ${Q}$ is $\Gamma_S(W)$-injective, then the implicit constant in the above estimate 
in Theorem \ref{th:l2_bound} is independent of $Q$.
}	
\end{rem}

\begin{proof}[Proof of Theorem \ref{th:l2_bound}]
Let again 
$$
\chi(g_\infty,g_S):=\chi_{E^{-1}}(g_\infty)\chi_{\cI_S}(g_S),\quad\quad \hbox{for $(g_\infty,g_S)\in G_\infty\times G_S$},
$$
denote the characteristic function 
of the subset $E^{-1}\times \cI_S$ of $G_\infty\times G_S$.
We introduce the subset 
$$
\Delta(W):=\Gamma_S\cap (G_\infty\times W)
$$
of $\Gamma_S$. Since $W$ is bi-invariant under the subgroup $U(W)$,
 it is clear that
$\Delta(W)$ is bi-invariant under the subgroup $\Gamma_S(W):=\Gamma_S\cap (G_\infty\times G_S\times U(W))$.
Moreover, since the set $W$ is union of finitely many right cosets of $U(W)$,
the set $\Delta(W)$ is a finite union of right cosets of $\Gamma_S(W)$.
For $g=(g_\infty,g_S)\in G_\infty\times G_S$, we define
$$
\phi(g):=\sum_{\delta\in \Delta(W)} \chi(g_\infty\delta^{-1},g_S\delta^{-1}).
$$ 
Since $\Delta(W)$ is $\Gamma_S(W)$-bi-invariant, this defines a measurable function on the space $Y_{S,W}=(G_\infty\times G_S)/\Gamma_S(W)$.
First, let us compute the integral of $\phi$:   
\begin{align}
\int_{Y_{S,W}} \phi \, d\mu_{S,W} &=
\int_{(G_\infty\times G_S)/\Gamma_S(W)} \left(\sum_{\delta\in \Delta(W)} \chi(g\delta^{-1})\right)\, d\mu_{S,W}(g)\nonumber \\
&=
\int_{(G_\infty \times G_S)/\Gamma_S(W)} \left(\sum_{\delta\in \Delta(W)/\Gamma_S(W)}\sum_{\gamma\in \Gamma_S(W)} \chi(g \gamma^{-1}\delta^{-1})\right)\, d\mu_{S,W}(g) \nonumber\\
&= \sum_{\delta\in \Delta(W)/\Gamma_S(W)}\int_{G_\infty \times G_S} \chi(g\delta^{-1})\, \frac{d(m_{\infty}\times m_{S})(g)}{|\Gamma_S/\Gamma_S(W)|} \nonumber\\
&= \frac{|\Delta(W)/\Gamma_S(W)|}{|\Gamma_S/\Gamma_S(W)|}  \int_{G_\infty \times G_S}  \chi\, d(m_{\infty}\times m_{S}) \nonumber \\
&= \frac{|\Delta(W)/\Gamma_S(W)|}{|\Gamma_S/\Gamma_S(W)|}  m_{\infty}(E^{-1})m_S(\cI_S)=
\frac{|\Delta(W)/\Gamma_S(W)|}{|\Gamma_S/\Gamma_S(W)|}  m_{\infty}(E). \nonumber
\end{align}
Since the projection of $\Gamma_S$ to $\cI^S$ is dense, the map $\gamma \to \gamma U(W)$, $\gamma\in \Gamma_S$, 
defines a bijection between the cosets 
$\Gamma_S/\Gamma_S(W)$ and $\cI^S/U(W)$.
In particular,
\begin{equation}\label{eq:v1}
|\Gamma_S/\Gamma_S(W)|=|\cI^S/U(W)|.
\end{equation}
Also, since $W$ is open, the projection of $\Delta(W)=\Gamma_S\cap (G_\infty\times G_S\times W)$ is dense in $W$, and we obtain that
\begin{equation}\label{eq:v2}
|\Delta(W)/\Gamma_S(W)|=|W/U(W)|.
\end{equation}
We recall  that the measure $m^S$ on $\cI^S$ is normalized so that $m^S(\cI^S)=1$. This implies that 
\begin{equation}\label{eq:v3}
m^S(W)=\frac{|W/U(W)|}{|\cI^S/U(W)|}.
\end{equation}
Hence, we conclude from the above computations that
\begin{equation}
\label{eq:l1}
\int_{Y_{S,W}} \phi \, d\mu_{S,W} =m_\infty(E)m^S(W).
\end{equation}
In particular, this shows that $\phi\in L^1(Y_{S,W})$.
In fact, we will show later in the proof that $\phi\in L^2(Y_{S,W})$.

Since $\phi$ in non-negative and integrable, it follows from the Fubini--Tonelli Theorem that for almost every $y\in Y_{S,W}$,
$$
\int_{a\in B_S(h)}\phi(ay)\, dm_{S}(a)<\infty.
$$

We shall show that the discrepancy $\mathcal{D}(R_S(h), E(\cdot)\times W)$ can be approximated by
such integrals. By the Fubini--Tonelli Theorem again, for $x\in G_\infty$ and $u\in \cI_S$,
\begin{align*}
\int_{a\in B_S(h)}\phi\big(a(x,u)\big)\, dm_{S}(a) &=\sum_{\delta\in \Delta(W)} \int_{a\in B_S(h)} \chi(x\delta^{-1}, au\delta^{-1})\, dm_{S}(a)\\
& = \sum_{\delta\in \Delta(W)} \int_{B_S(h)} \chi_{E^{-1}}(x\delta^{-1})\chi_{\cI_S} (au\delta^{-1})  \, dm_{S}(a)\\
&=
\sum_{\delta\in \Delta(W)\cap (Ex\times G_S)} m_{S}\big(\cI_S \delta u^{-1}\cap B_S(h)\big).
\end{align*}
We recall that the sets $B_S(h):=\{b\in G_S:\, \h(b)\le h\}$ are $\cI_S$-bi-invariant.  Therefore, 
if $\delta\in B_S(h)$, we have $\cI_S\delta u^{-1}\subset B_S(h)$, and
if $\delta\notin B_S(h)$, we have $\cI_S\delta u^{-1}\cap B_S(h)=\emptyset$.
Hence, since $m_S(\cI_S)=1$, it follows that for every $u\in \cI_S$,
\begin{align}\label{eq:form}
\big|R_S(h)\cap (E(x)\times W)\big|=\int_{a\in B_S(h)}\phi\big(a(x,u)\big)\, dm_{S}(a).
\end{align}
Combining \eqref{eq:l1} and \eqref{eq:form}, we derive that
\begin{align*}
\mathcal{D}(R_S(h), E(x)\times W) &=
\left| \frac{1}{m_S(B_S(h))}\int_{a\in B_S(h)}\phi\big(a(x,e)\big)\, dm_{S}(a)-\int_{Y_{S,W}} \phi\, d\mu_{S,W}  \right|\\
&=\left| \pi_{{S,W}}(\beta_h) \phi (x,e)-\int_{Y_{S,W}} \phi\, d\mu_{S,W}  \right|.
\end{align*}
This representation allows us to apply the estimate \eqref{eq:nnorm}.
Let $Q$ be a bounded measurable subset of $G_\infty$.
Since the integral in \eqref{eq:form} is independent of $u\in \cI_S$,
we obtain that
\begin{align*}
&\int_Q \left| \int_{a\in B_S(h)}\phi\big(a(x,e)\big)\, dm_{S}(a)-m_S(B_S(h))\int_{Y_{S,W}} \phi\, d\mu_{S,W}  \right|^2 \, dm_\infty(x)\\
= & \int_{Q\times \cI_S} \left| \int_{a\in B_S(h)}\phi\big(a(x,u)\big)\, dm_{S}(a)-m_S(B_S(h))\int_{Y_{S,W}} \phi\, d\mu_{S,W}  \right|^2 \, dm_\infty(x)dm_S(u).
\end{align*}
The set $Q\times \cI_S$ projects to a measurable subset $(Q\times \cI_S)\Gamma_S(W)$ in $Y_{S,W}$. Since this is a bounded subset of $G_\infty\times G_S$, the fibers of this map have cardinalities bounded by a uniform constant $N_Q$, which is independent of $W$ because $\Gamma_S(W)\subset \Gamma_S$. Hence, the last integral is bounded by
\begin{align*}
& N_Q\, \int_{(Q\times \cI_S)/\Gamma_S(W)} \left| \int_{a\in B_S(h)}\phi(ay)\, dm_{S}(a)-m_S(B_S(h))\int_{Y_{S,W}} \phi\, d\mu_{S,W}  \right|^2 \, d\mu_{S,W}(y)\\
\le \,& N_Q\, \int_{Y_{S,W}} \left| \int_{a\in B_S(h)}\phi(ay)\, dm_{S}(a)-m_S(B_S(h))\int_{Y_{S,W}} \phi\, d\mu_{S,W}  \right|^2 \, d\mu_{S,W}(y)\\
=& \, N_Q\, m_S(B_S(h))^2
\left\| \pi_{{S,W}}(\beta_h) \phi -\int_{Y_{S,W}} \phi \, d\mu_{S,W} \right\|^2_{L^2(Y_{S,W})}.
\end{align*}

Thus, it follows from \eqref{eq:nnorm} that for every $\eta > 0$ 
\begin{equation}\label{eq:prev}
\Big\|\mathcal{D}(R_S(h), E(\cdot)\times W)\Big\|_{L^2(Q)}\ll_{S,Q,\eta} 
\, m_S(B_S(h))^{-\mathfrak{k}_S+\eta}\, \|\phi\|_{L^2(Y_{S,W})}.
\end{equation}
The function $\phi$ is indeed square-integrable, and the $L^2$-norm $\|\phi\|_{L^2(Y_{S,W})}^2$ is computed as follows. 
\begin{align*}
 & \int_{(G_\infty\times G_S)/\Gamma_S(W)} \left(\sum_{\delta_1,\delta_2\in \Delta(W)} \chi(g\delta_1^{-1})\chi(g\delta_2^{-1}) \right)\, d\mu_{S,W}(g)\\
=&
\int_{(G_\infty\times G_S)/\Gamma_S(W)} \left(
\sum_{\delta_1,\delta_2\in \Delta(W)/\Gamma_S(W)}\sum_{\gamma_1,\gamma_2\in \Gamma_S(W)} \chi(g\gamma_1^{-1}\delta_1^{-1})\chi(g\gamma_2^{-1}\delta_2^{-1}) \right)\, d\mu_{S,W}(g)\\
= &
\int_{(G_\infty\times G_S)/\Gamma_S(W)} \left(
\sum_{\delta_1,\delta_2\in \Delta(W)/\Gamma_S(W)}\sum_{\sigma,\gamma\in \Gamma_S(W)} \chi(g\sigma\delta_1^{-1})\chi(g\sigma\gamma^{-1}\delta_2^{-1}) \right)\, d\mu_{S,W}(g)\\
= &
\int_{G_\infty\times G_S} \left(
\sum_{\delta_1,\delta_2\in \Delta(W)/\Gamma_S(W)}\sum_{\gamma\in \Gamma_S(W)} \chi(g\delta_1^{-1})\chi(g\gamma^{-1}\delta_2^{-1}) \right)\, \frac{d(m_\infty\times m_S)(g)}{|\Gamma_S:\Gamma_S(W)|}\\
= &
\int_{G_\infty\times G_S} \left(
\sum_{\gamma\in \Delta(W)/\Gamma_S(W)}\sum_{\delta\in \Delta(W)} \chi(g\gamma^{-1})\chi(g\delta^{-1}) \right)\, \frac{d(m_\infty\times m_S)(g)}{|\Gamma_S:\Gamma_S(W)|}.
\end{align*}
We observe that for any $\gamma,\delta\in \Delta(W)$,
\begin{align*}
\int_{G_\infty\times G_S} \chi(g\gamma^{-1})\chi(g\delta^{-1})\, d(m_\infty\times m_S)(g)&=
m_\infty(E^{-1}\gamma\cap E^{-1} \delta)m_S(\cI_S\gamma\cap \cI_S \delta) \\
&\le m_\infty(E)m_S(\cI_S)=m_\infty(E). 
\end{align*}
Moreover this integral is zero unless $m_\infty(E^{-1}\cap E^{-1}\delta\gamma^{-1})> 0$
and $\delta\gamma^{-1}\in \cI_S$. But since $\gamma, \delta \in \Delta(W)$, we also have $\delta\gamma^{-1}\in \cI^S$, and so in this case it follows that $\delta\gamma^{-1}\in {\sf G}(O)$.
Recalling the definition of $\mathcal{N}(E)$, this implies that for fixed $\gamma\in\Delta(W)$, 
$$
\left|\left\{\delta\in \Delta(W): \int_{G_\infty\times G_S} \chi(g\gamma^{-1})\chi(g\delta^{-1})\, d(m_\infty\times m_S)(g)\ne 0\right\}\right|\le \mathcal{N}(E).
$$
Hence, applying the Fubini--Tonelli Theorem once more, we conclude that
$$
\|\phi\|_{L^2(Y_{S,W})}^2\le \mathcal{N}(E)m_\infty(E)\frac{|\Delta(W)/\Gamma_S(W)|}{|\Gamma_S:\Gamma_S(W)|}=\mathcal{N}(E)m_\infty(E)m^S(W).
$$

Here we have used \eqref{eq:v1}--\eqref{eq:v3} in the last step.
Combining this estimate with \eqref{eq:prev}, this completes the proof of Theorem \ref{th:l2_bound}.
\end{proof}

\subsection{Almost sure discrepancy estimates}

We now turn to establish an almost sure bound for the discrepancy. 
We will use the mean-square  bound of Theorem \ref{th:l2_bound} 
and the Borel-Cantelli lemma, and the main difficulty here will be to establish an estimate which holds for all 
$h$ on a fixed set of full measure. For this argument, we need the following elementary lemma:

\begin{lem}\label{lem:elm}
For $\theta\ge 1$ and $y_1,\ldots, y_n\ge 0$,
$y_1^\theta+\cdots + y_n^\theta\le (y_1+\cdots + y_n)^\theta.$
\end{lem}

\begin{proof}
The proof proceeds by induction on $n$. We consider the function
$$
f(x):=x^\theta+y_2^\theta\cdots + y_n^\theta- (x+y_2+\cdots + y_n)^\theta.
$$
By the inductive assumption, $f(0)\le 0$. Using that $\theta\ge 1$, one checks that $f'(x)\le 0$ for $x\ge 0$. This implies the claim. 	
\end{proof}

\begin{thm}[almost sure discrepancy bound]\label{almost-sure-disc}
With notation as in Theorem~\ref{th:l2_bound}, for every $0< \mathfrak{k}< \mathfrak{k}_S({\sf G})$  for almost all $x\in G_\infty$, and for every $\eta > 0$ 
$$
\mathcal{D}\big(R_S(h), E(x)\times W\big)\ll_{S,E,W,x,\mathfrak{k}, \eta}  \big(\log v_S(h)\big)^{3/2+\eta} v(h)^{-\mathfrak{k}}.
$$	
\end{thm}

\begin{proof}
For $a, b\in \mathbb{Z}_{\ge 0}$, we consider the intervals
$$
I_{a,b}:=\big\{h\ge 1:\, 2^b a <v_S(h)\le 2^b (a+1) \big\}.
$$
Note that for fixed $b$, they define a partition of $[1,\infty)$.

For $s\in\mathbb{Z}_{\ge 0}$, we set
$$
\mathcal{M}_{s,b}:=\{I_{a,b}:\, v_S(I_{a,b})\subset (0,2^s] \}\quad\hbox{and}\quad
\mathcal{M}_{s}:=\sqcup_{b\le s}\,  \mathcal{M}_{s,b}.
$$
We observe that 
$$
\bigcup_{I\in \mathcal{M}_{s}} I= [1,h_s]\quad \hbox{with $2^{s-1}< v_S(h_s)\le 2^s$}. 
$$
For $I\subset [1,\infty)$, we set 
$$
R_S(I):=\{\gamma\in \Gamma_S:\, \h(\gamma)\in I \}\quad\hbox{and}\quad
B_S(I):=\{g\in G_S:\, \h(g)\in I \}.
$$
The argument of the proof of Theorem \ref{th:l2_bound} gives the following mean bound:
for bounded measurable subsets $Q\subset G_\infty$, 
\begin{align*}
\int_Q \Big| \big|R_S(I)\cap (E(x)\times W)\big|-m_\infty(E)m^S(W)m_S(B_S(I))\Big|^2 \, dm_\infty(x)
\ll  \,  m_S(B_S(I))^{2-2\mathfrak{k}}.
\end{align*}
The implicit constant here and the following computations depends on $S,E,W,Q$ and the difference $\mathfrak{k}_S-\mathfrak{k}$.
Using Lemma~\ref{lem:elm}, 
$$
\sum_{I\in\mathcal{M}_{s,b}} m_S(B_S(I))^{2-2\mathfrak{k}}\le 
\left(\sum_{I\in\mathcal{M}_{s,b}} m_S(B_S(I))\right)^{2-2\mathfrak{k}}
\le m_S(B_S(h_s))^{2-2\mathfrak{k}} \le 2^{s(2-2\mathfrak{k})}.
$$	
Hence,
\begin{align}
\sum_{I\in\mathcal{M}_{s}} \int_Q \Big|R_S(I)\cap (E(x)\times W)|  -m_\infty(E)m^S(W)m_S(B_S(I)) \Big|^2 \,  dm_\infty(x) 
\ll  \;  s 2^{s(2-2\mathfrak{k})}.\label{eq:sum_bound}
\end{align}

Let $\eta>0$. We consider the sets $\Upsilon_s$ consisting of $x\in Q$ 
such that 
\begin{align*}
\sum_{I\in\mathcal{M}_{s}} \Big|R_S(I)\cap (E(x)\times W)|-m_\infty(E)m^S(W)m_S(B_S(I))&\Big|^2 
\ge\,  s^{2+\eta} 2^{s(2-2\mathfrak{k})}.
\end{align*}
We deduce from \eqref{eq:sum_bound} that 
$$
m_\infty(\Upsilon_s)\ll s^{-1-\eta}.
$$
It follows from the Borel--Cantelli Lemma that the $\limsup$ of the sets $\Upsilon_s$ has measure zero. Hence, for almost all $x\in Q$ and all $s\ge s_0(x)$, we have the bound
\begin{align}
\sum_{I\in\mathcal{M}_{s}} \Big|\big|R_S(I)\cap (E(x)\times W)\big|-m_\infty(E)m^S(W)m_S(B_S(I))\Big|^2 
\ll\,  s^{2+\eta} 2^{s(2-2\mathfrak{k})}.\label{eq:sum0}
\end{align}
We shall use this bound to prove the theorem.

We first consider the case when $h$ is an end point of one of the intervals $I_{a,b}$.
We choose the parameter $s$ such that $2^{s-1}< v_S(h)\le 2^s$.
Using the binary representation, the interval $[1,h]$ can be written 
as a disjoint intervals of at most $s$ intervals $I_i$ from $\mathcal{M}_s$.
Then using \eqref{eq:sum0} and the Cauchy--Schwartz inequality, we deduce that 
for almost all $x\in Q$ and $h\ge h_0(x)$, and for every $\eta > 0$ 
\begin{align*}
\Big|\big|R_S(h)\cap (E(x)\times W)\big|&-m_\infty(E)m^S(W)m_S(B_S(h))\Big|\\
&\le \sum_i \Big|\big|R_S(I_i)\cap (E(x)\times W)\big|-m_\infty(E)m^S(W)m_S(B_S(I_i))\Big|\\
&\ll_\eta \, s^{3/2+\eta/2} 2^{s(1-\mathfrak{k})}\ll_\eta \, (\log v_S(h))^{3/2+\eta/2} v_S(h)^{1-\mathfrak{k}}.
\end{align*}
For a general $h$, we observe that there exists $h_1$ and $h_2$ as above such that
$$
h_1\le h\le h_2\quad\hbox{and}\quad v_S(h_2)-v_S(h_1)\le 1.
$$
Then 
$$
\big|R_S(h_1)\cap (E(x)\times W)\big|\le \big|R_S(h)\cap (E(x)\times W)\big|\le \big|R_S(h_2)\cap (E(x)\times W)\big|,
$$
and 
$$
v_S(h_2)\le v_S(h)+1\quad\hbox{and}\quad v_S(h_1)\ge v_S(h)-1.
$$
This provides upper and lower bounds on $|R_S(h)\cap (E(x)\times W)|$
that imply that for almost all $x\in Q$ and all sufficiently large $h$,
\begin{align*}
\big|R_S(h)\cap (E(x)\times W)\big|=\, m_\infty(E)m^S(W)v_S(h)
 + O\Big( (\log v_S(h))^{3/2+\eta/2} v_S(h)^{1-\mathfrak{k}}\Big).
\end{align*}
Clearly, this estimate also holds for all $h$ with an implicit constant depending on $x$ and $\eta$.
Finally, exhausting $G_\infty$ by a countable union
of bounded measurable sets, we deduce that 
this estimate holds for almost all $x\in G_\infty$,
which implies Theorem \ref{almost-sure-disc}. 

\end{proof}

\subsection{Uniform discrepancy estimates}
We now turn to establish a uniform pointwise bound on the discrepancy of rational points in the collection of sets $E(x)\times W$ with $E(x):=Ex\subset G_\infty$ and $W\subset \cI^S$. 
In our discussion the set $E$ is fixed, the congruence condition $W$ is arbitrary, and we consider {\it  all points }$x\in G_\infty$. 
This requires imposing a regularity condition on the set $E$, namely
the right-stability property \eqref{eq:well}.  The bound will now involve the dimension $\mathfrak{d}:=\dim_\RR (G_\infty)$, whereas the almost sure pointwise bound in Theorem \ref{almost-sure-disc} did not. 

 We recall that by Lemma \ref{l:ball}(a),
\begin{equation} \label{eq:vv2}\tag{V}
\epsilon^{\mathfrak{d}}\ll m_\infty(B(e,\epsilon)) \ll \epsilon^{\mathfrak{d}}\quad\hbox{for all $\epsilon\in (0,\epsilon_0)$.}
\end{equation}

\begin{thm}[uniform discrepancy bound]\label{th:wr_general}

With notation as in Theorem \ref{th:l2_bound}, 
let $E$ be a right-stable finite-measure subset of $G_\infty$ satisfying $\mathcal{N}(E_{\epsilon_0}^+)<\infty$, and $W$ compact open subset of $\cI^S$. For every 
$0 < \mathfrak{k}= \mathfrak{k}_S({\sf G})-\eta < \mathfrak{k}_S({\sf G})$ 
and for $x\in G_\infty$ the following pointwise bound for the discrepancy holds :

$$
	\mathcal{D}(R_S(h), E(x)\times W) \ll_{S,E,x,\eta} m^S(W)^{(\mathfrak{d}+1)/(\mathfrak{d}+2)}\, v_S(h)^{-2\mathfrak{k}/(\mathfrak{d}+2)}\,.
	$$
provided that 
$m^S(W)\gg_\eta \, v_S(h)^{-2\mathfrak{k}}.$

This condition is also equivalent  to the condition that the error term in the asymptotic expansion the solution counting function :
\begin{align*}
\big|R_S(h)\cap (E(x)\times W)\big|=&\,m_\infty(E)m^S(W)v_S(h)\\
&+ O_{S,E,x,\eta}\Big( m^S(W)^{(\mathfrak{d}+1)/(\mathfrak{d}+2)}\, v_S(h)^{1-2\mathfrak{k}/(d+2)}\Big).
\end{align*}
is bounded by the main term.

Explicity, if the volume growth satisfies $v_S(h)\gg_{S,\eta^\prime} h^{\mathfrak{a}}$, with $0< \mathfrak{a}=\mathfrak{a}_S({\sf G})-\eta^\prime$, then the estimate holds provided the height $h$ satisfied $h\gg_{\eta^\prime, \eta}  (1/m^S(W))^{1/2\mathfrak{a}\mathfrak{k}}.$
Moreover, the above  estimate is uniform for $x$ ranging in compact subsets of $G_\infty$.

\end{thm}

In the proof, we use the estimates \eqref{eq:nnorm}, \eqref{eq:well}, and \eqref{eq:vv2}.

\begin{proof}[Proof of Theorem \ref{th:wr_general}]
The starting point of our argument is the $L^2$-bound established in Theorem \ref{th:l2_bound}. 
Let $Q$ be a bounded measurable subset of $G_\infty$,
and define the open set $Q':=B(e,\epsilon_0)Q.$
According to Theorem \ref{th:l2_bound}, the following bound holds:  
for any measurable subset $F$ of $G_\infty$ with finite positive measure, and any given $\eta > 0$:
$$
\int_{Q'}\mathcal{D}\big(R_S(h), F(y)\times W\big)^2\, dm_\infty(y)\le C\, \mathcal{N}(F) m_\infty(F) m^S(W) v_S(h)^{-2\mathfrak{k}}
$$
for some $C=C_{Q, S, \epsilon_0, \eta}>0$.
This implies that for every $\delta>0$,
\begin{align}
m_\infty\Big(\big\{y\in Q':\, \mathcal{D}(R_S(h), & F(y)\times W)  >\delta \big\}\Big) \nonumber \\
&\le\, C\, \delta^{-2}\mathcal{N}(F)m_\infty(F) m^S(W) v_S(h)^{-2\mathfrak{k}}. \label{eq:b_Q}
\end{align}
We introduce a parameter $\epsilon\in (0,\epsilon_0)$ which we assume to satisfy
\begin{equation}\label{eq:delta_00}
m_\infty\big(B(e,\epsilon)\big)>C\, \delta^{-2}\mathcal{N}(E^+_\epsilon)m_\infty(E_\epsilon^+) m^S(W) v_S(h)^{-2\mathfrak{k}}. 
\end{equation}
Here $\epsilon_0$ is determined by conditions \eqref{eq:well} and \eqref{eq:vv2}.
We note that for $x\in Q$, we have $B(e,\epsilon)x\subset Q'$.
Hence, it follows from \eqref{eq:b_Q}
applied to $F=E_\epsilon^+$ that for every $x\in Q$,
there exists $g\in B(e,\epsilon)$ such that 
\begin{equation}\label{eq:up}
\mathcal{D}\big(R_S(h), E_\epsilon^+(gx)\times W\big)\le\delta.
\end{equation}

We observe when $g\in B(e,\epsilon)=B(e,\epsilon)^{-1}$, it follows from the definition of $E_\epsilon^-$
and $E_\epsilon^+$   (cf.~\eqref{eq:well}) that
\begin{equation}\label{eq:www}
E_\epsilon^-(gx)\subset E(x)\subset E^+_\epsilon(gx).
\end{equation}
Hence, taking $g\in B(e,\epsilon)$ as in \eqref{eq:up}, we obtain that
\begin{align*}
\big|R_S(h)\cap \big( E(x)\times W\big) \big|&\le \big|R_S(h)\cap \big(E_\epsilon^+(gx)\times W\big)\big|\\
&\le \Big(
m_\infty(E_\epsilon^+(gx)) m^S(W)+\mathcal{D}\big(R_S(h), E_\epsilon^+(gx)\times W\big)\Big)v_S(h)\\
&\le \big(m_\infty(E_\epsilon^+) m^S(W)+\delta\big)v_S(h).
\end{align*}
Since the set $E$ is assumed to be right stable,
we deduce that
\begin{align*}
\big|R_S(h)\cap \big( E(x)\times W\big) \big|\le
\left(
m_\infty(E) m^S(W)+O(\epsilon\, m^S(W))+\delta\right)v_S(h)
\end{align*}
for all $\epsilon\in (0,\epsilon_0)$ satisfying \eqref{eq:delta_00}.

We also apply a similar argument to deduce a lower bound.
Since it follows from \eqref{eq:delta_00} that also
$$
m_\infty(B(e,\epsilon))>C\, \delta^{-2}\mathcal{N}(E^-_\epsilon)m_\infty(E_\epsilon^-) m^S(W) v(h)^{-2\mathfrak{k}}, 
$$
we deduce from \eqref{eq:b_Q} with $F=E_\epsilon^-$  that for every $x\in Q$,
there exists $g\in B(e,\epsilon)$ such that
$$
\mathcal{D}\big(R_S(h), E_\epsilon^-(gx)\times W\big)\le\delta.
$$
Hence, we deduce as above that
\begin{align*}
\big|R_S(h)\cap \big(E(x)\times W\big)\big|
&\ge \big|R_S(h)\cap \big(E_\epsilon^-(gx)\times W\big)\big|\\
&\ge \Big(m_\infty(E_\epsilon^-(gx)) m^S(W)-\mathcal{D}\big(R_S(h), E_\epsilon^-(gx)\times W\big)
\Big)v_S(h)\\
&\ge
\left(m_\infty(E) m^S(W) -O\big(\epsilon\,m^S(W)
\big)+\delta\right)v_S(h).
\end{align*}

Combining the above estimates on $\big|R_S(h)\cap \big(E(x)\times W\big)\big|$, we deduce that for every $x\in Q$,
$$
\mathcal{D}\big(R_S(h), E(x)\times W\big)\le \delta+O\big(\epsilon\,
m^S(W)\big)
$$
provided that \eqref{eq:delta_00} is satisfied. 
To arrange \eqref{eq:delta_00}, it is sufficient to pick $\delta$ of the form
$$
\delta=c\,m_\infty(B(e,\epsilon))^{-1/2}\mathcal{N}(E_\epsilon^+)^{1/2}m_\infty(E_\epsilon^+)^{1/2} m^S(W)^{1/2} v_S(h)^{-\mathfrak{k}}
$$
with sufficiently large $c=c_{Q,S, \epsilon_0, \eta}>0$. Then we deduce using \eqref{eq:vv2} that for all $x\in Q$,
\begin{equation}\label{eq:bd}
\mathcal{D}\big(R_S(h), E(x)\times W\big)\ll_{E,Q} \epsilon^{-\mathfrak{d}/2} m^S(W)^{1/2} v_S(h)^{-\mathfrak{k}}+
\epsilon\,m^S(W).
\end{equation}
We balance the two summands in the estimate and choose the parameter $\epsilon$ as 
$$
\epsilon=\Big( m^S(W)^{-1/2} v_S(h)^{-\mathfrak{k}}\Big)^{1/(\mathfrak{d}/2+1)}\,.
$$
Therefore, when $h$ is sufficiently large, we have $\epsilon\in (0,\epsilon_0)$ provided that 
$$
v_S(h)^{2\mathfrak{k}}m^S(W) > \epsilon_0^{-(\mathfrak{d}+2)},  
\text{  namely }  v_S(h)^{2\mathfrak{k}}\gg m^S(W)^{-1}\,.
$$
the latter condition is equivalent to 
$$
m^S(W)\gg m^S(W)^{(\mathfrak{d}+1)/(\mathfrak{d}+2)}\, v_S(h)^{-2\mathfrak{k}/(\mathfrak{d}+2)}, 
$$
and substituting this expression in the bound just established for the discrepancy  we deduce that for every $x\in Q$,
	$$
\mathcal{D}\big(R_S(h), E(x)\times W\big) \ll_{S,E,Q,\eta} m^S(W)^{(\mathfrak{d}+1)/(\mathfrak{d}+2)}\, v_S(h)^{-2\mathfrak{k}/(\mathfrak{d}+2)}.
$$
This completes the proof of Theorem \ref{th:wr_general}.
\end{proof}

\subsection{Uniform discrepancy at arbitrarily small scales}
We now turn to establish a bound on the discrepancy which is uniform over a family of balls $B(x, \ell)$ in $G_\infty$ of arbitrarily small radius, for {\it every} $x\in G_\infty$. Let $\ell>0$, and $W$ be a compact open subset of $\cI^S$.
We will consider the family of subsets $\Omega$ given by 
$
B(x,\ell)\times W\subset  G_\infty\times \cI^S, 
$
and our goal is to prove an explicit pointwise bound
for the number of $\Gamma_S$-points contained in them.

\begin{thm}[pointwise discrepancy bound at small scales]\label{th:balls_general}
With notation as in Theorem \ref{th:l2_bound},  let $W$ be a compact open subset of $\cI^S$. Fix  
$0 < \mathfrak{k}= \mathfrak{k}_S({\sf G})-\eta < \mathfrak{k}_S({\sf G})$,  
let $x\in G_\infty$ and $\ell\in (0,\ell_0)$ (for suitable $\ell_0 >0$ independent of $x$), and set 
$$
\mathcal{E}_{\ell,W}(h):= m_\infty(B(e,\ell))^{\mathfrak{d}/(\mathfrak{d}+2)}m^S(W)^{(\mathfrak{d}+1)/(\mathfrak{d}+2)}v_S(h)^{1-2\mathfrak{k}/(\mathfrak{d}+2)}.
$$
Then 
\begin{equation}\label{eq:main}
\big|R_S(h)\cap B(x,\ell)\times W\big|= m_\infty(B(e,\ell))m^S(W)v_S(h)+O_{S,x}\Big (\mathcal{E}_{\ell,W}(h)\Big)
\end{equation}
provided that 
\begin{equation}\label{eq:condition} 
m_\infty(B(e,\ell))^2 m^S(W)\ge v_S(h)^{-2\mathfrak{k}}\,.
\end{equation} 
This condition is equivalent to the condition that the error term in the foregoing asymptotic formula is bounded by the main term. 
%
Explicitly, if $v_S(h) \gg_{\eta^\prime} h^{\mathfrak{a}}$ with $0< \mathfrak{a}={\mathfrak{a}_S({\sf G})-\eta^\prime}$ for $\eta^\prime > 0$, 
then the result holds at any scale $0< \ell < \ell_0$, provided $h$ satisfies 
$$h\gg_{S,\eta}\ell^{-\mathfrak{d}/\mathfrak{k}\mathfrak{a}} m^S(W)^{-(\mathfrak{d}+2)/2\mathfrak{k}\mathfrak{a}}\,.$$
Moreover, this estimate is uniform for $x$ ranging in compact subsets of $G_\infty$.
\end{thm}

\begin{proof} 
We adopt the argument from the proof of Theorem \ref{th:wr_general} with some modifications.
We note that there exists  $r_0 > 0$ such that when $r\in (0,r_0)$,
$$
B(e,r)\cap \gamma B(e,r)=\emptyset\quad\;\hbox{ for $\gamma\in {\sf G}(O)\backslash \{e\}$,}
$$
so that $\mathcal{N}(B(e,r))=1$.

We fix a bounded measurable subset $Q$ of $G_\infty$
and define the open set $Q':=B(e,r_0)Q$.
Applying Theorem \ref{th:l2_bound}
to the balls $B(y,\ell)=B(e,\ell)y$, we deduce that for every $r\in (0,r_0)$ and $\delta>0$,	and some $C=C_{S,Q, r_0}$ 
\begin{align}
m_\infty\Big(\big\{y\in Q':\, \mathcal{D}\big(R_S(h), B(x,\ell)\times W\big)&>\delta \big\}\Big) \nonumber \\
\le\, & C\, \delta^{-2}m_\infty(B(e,r)) m^S(W) v_S(h)^{-2\mathfrak{k}}. \label{eq:m_inf}
\end{align}
We introduce a positive parameter $\epsilon$ (to be chosen later) satisfying 
\begin{equation}\label{eq:delta_0}
\epsilon<\ell,
\end{equation} 
as well as 
\begin{equation}\label{eq:delta}
m_\infty\big(B(e,\epsilon)\big)>C\, \delta^{-2}m_\infty(B(e,\ell+\epsilon)) m^S(W) v_S(h)^{-2\mathfrak{k}}. 
\end{equation}
We choose $\ell_0$ so that $0 < 2\ell_0< r_0$
and also so that Lemma \ref{l:ball} is applicable for $r < r_0$.
Then it follows from \eqref{eq:m_inf} that for every $x\in Q$,
there exists $g\in B(e,\epsilon)$ such that 
$$
\mathcal{D}\big(R_S(h), B(gx,\ell+\epsilon)\times W\big)\le\delta.
$$
It follows from the triangle inequality that
$$
B(gx,\ell-\epsilon)\times W\subset B(x,\ell)\times W\subset B(gx,\ell+\epsilon)\times W,
$$
and so,
\begin{align*}
\big|R_S(h)\cap (B(x,\ell)\times W)\big|&\le \big| R_S(h)\cap
(B(gx,\ell+\epsilon)\times W)\big|\\
&\le \Big(
m_\infty(B(e,\ell+\epsilon)) m^S(W)+\mathcal{D}\big(R_S(h), 
B(gx,\ell+\epsilon)\times W\big)\Big)v_S(h)\\
&\le \big(
m_\infty(B(e,\ell+\epsilon)) m^S(W)+\delta\big)v_S(h).
\end{align*}
Since it follows from Lemma \ref{l:ball}(b) that  
\begin{equation}\label{eq:ball}
m_\infty(B(e,\ell+\epsilon))\le \left(1+c\,\frac{\epsilon}{\ell}\right)m_\infty(B(e,\ell)), 
\end{equation}
we deduce that
\begin{align*}
\big|R_S(h)\cap (B(x,\ell)\times W)\big|\le
\left( 
m_\infty(B(e,\ell)) m^S(W)+\delta+c\,\frac{\epsilon}{\ell} m_\infty(B(e,\ell))m^S(W)\right)v_S(h).
\end{align*}
A similar argument also gives the  lower bound 
\begin{align*}
\big|R_S(h)\cap (B(x,\ell)\times W)\big|\ge
\left(
m_\infty(B(e,\ell)) m^S(W)-\delta-c\,\frac{\epsilon}{\ell} m_\infty(B(e,\ell))m^S(W)\right)v_S(h).
\end{align*}
It follows that the discrepancy 
$$
\mathcal{D}(R_S(h), B(x,\ell)\times W)=\left| \frac{\big|R_S(h)\cap (B(x,\ell)\times W)\big|}{v_S(h)} - m_\infty(B(e,\ell)) m^S(W) \right|.
$$
satisfies
$$
\mathcal{D}(R_S(h), B(x,\ell)\times W)\le \delta+
c\,\frac{\epsilon}{\ell} m_\infty(B(e,\ell))m^S(W).
$$
The parameter $\delta$ have to satisfy \eqref{eq:delta}, 
and since $m_\infty(B(e,\ell+\epsilon))\le (1+c)m_\infty(B(e,\ell))$ (using \eqref{eq:ball} and $\epsilon < \ell$) it is sufficient to pick $\delta$ of the form
$$
\delta=b_0\,m_\infty(B(e,\epsilon))^{-1/2} m_\infty(B(e,\ell))^{1/2} m^S(W)^{1/2} v_S(h)^{-\mathfrak{k}}
$$
with sufficiently large $b_0>0$. Then using Lemma \ref{l:ball}(a),
we deduce the bound
\begin{equation}\label{eq:bd}
 \mathcal{D}(R_S(h), B(x,\ell)\times W)\ll \epsilon^{-\mathfrak{d}/2} m_\infty(B(e,\ell))^{1/2} m^S(W)^{1/2} v_S(h)^{-\mathfrak{k}}+
\frac{\epsilon}{\ell} m_\infty(B(e,\ell))m^S(W).
\end{equation}
The two summands in the foregoing bound are balanced precisely when $\epsilon$ is a constant multiple of 
\begin{equation}\label{eq:epsilon}
\ell^{2/(\mathfrak{d}+2)}m_\infty(B(e,l))^{-1/(\mathfrak{d}+2)}m^S(W)^{-1/(\mathfrak{d}+2)}v_S(h)^{-2\mathfrak{k}/(\mathfrak{d}+2)}\,,
\end{equation}
and then the estimate that arises in (\ref{eq:bd}) is give by a constant multiple of 
$$
\ell^{-\mathfrak{d}/(\mathfrak{d}+2)}m_\infty(B(e,l)^{(\mathfrak{d}+1)/(\mathfrak{d}+2)}m^S(W)^{(\mathfrak{d}+1)/(\mathfrak{d}+2) }v_S(h)^{-2\mathfrak{k}/(\mathfrak{d}+2)}\,.
$$
By Lemma \ref{l:ball}(a) 
$\ell \ll m_\infty(B(e,l)^{1/\mathfrak{d}}\ll \ell$, and so let us define the error term in (\ref{eq:main}) by 
\begin{equation}\label{eq:err}
\mathcal{E}_{\ell,W}(h):= m_\infty(B(e,\ell))^{\mathfrak{d}/(\mathfrak{d}+2)}m^S(W)^{(\mathfrak{d}+1)/(\mathfrak{d}+2)}v_S(h)^{1-2\mathfrak{k}/(\mathfrak{d}+2)}\,.
\end{equation}

We now claim that under the condition (\ref{eq:condition}) stated in Theorem \ref{th:balls_general} namely $m_\infty(B(e,\ell))^2 m^S(W)\ge v_S(h)^{-2\mathfrak{k}}$, the conclusion of Theorem \ref{th:balls_general} holds: 
\begin{equation}\label{eq:mmm}
\big|R_S(h)\cap (B(e,\ell)\times W)\big|= m_\infty(B(e,\ell))m^S(W)v_S(h)+O_Q\Big (\mathcal{E}_{\ell,W}(h)\Big).
\end{equation}
Indeed,  (\ref{eq:condition}) is equivalent to the condition that the error term $\mathcal{E}_{\ell,W}(h)$ defined in (\ref{eq:err}) is bounded by the main term in (\ref{eq:mmm}), namely :
\begin{equation}\label{eq:eee}
 m_\infty(B(e,\ell))m^S(W)v_S(h)\ge \mathcal{E}_{\ell,W}(h)\,.
\end{equation}
Furthermore, these two conditions are equivalent to the condition 
\begin{equation}\label{eq:eee_0}
 m_\infty(B(e,\ell))^{-2/(\mathfrak{d}+2)} m^S(W)^{-1/(\mathfrak{d}+2)}v_S(h)^{-2\mathfrak{k}/(\mathfrak{d}+2)}\le 1.
\end{equation}
To establish (\ref{eq:mmm}) we apply \eqref{eq:bd} with the choice of the parameter $\epsilon$ being given by (\ref{eq:epsilon}):
\begin{align*}
\epsilon&=c_0'\, \ell^{2/(\mathfrak{d}+2)} m_\infty(B(e,\ell))^{-1/(\mathfrak{d}+2)} m^S(W)^{-1/)(\mathfrak{d}+2)} v_S(h)^{-2\mathfrak{k}/(\mathfrak{d}+2)}.
\end{align*}
where $c_0'>0$ is a sufficiently small positive constant so that Lemma \ref{l:ball}(b) holds. 
It follows from Lemma \ref{l:ball}(a) and \eqref{eq:eee_0} that
\begin{align*}
\epsilon
\ll c_0' \,\ell^{-(\mathfrak{d}-2)/(\mathfrak{d}+2)} m^S(W)^{-1/(\mathfrak{d}+2)} v_S(h)^{-2\mathfrak{k}/(\mathfrak{d}+2)}\ll c_0'\,\ell.
\end{align*}
Hence, with sufficiently small $c_0'$, condition \eqref{eq:delta_0} holds, 
and so the above estimates apply.
Now the bound \eqref{eq:mmm} follows from \eqref{eq:bd} and this completes the proof of Theorem \ref{th:balls_general}.
\end{proof}

\section{Discrepancy bounds for general groups}\label{sec:gen}
In the present section, our goal is to extend the results of the previous sections to the general case, where 
$\sf G$ is {\it any} connected $K$-simple algebraic group defined over a   number field $K$. As above, we fix a subset $S$ of finite places of $K$ such that $\sf G$ is isotropic over $S$, consider the $S$-arithmetic group $\Gamma_S:={\sf G}(O_S)$, and aim to analyze the distribution of $\Gamma_S$ embedded in $G_\infty$. 

More generally, the method that we develop can be used to analyze the distribution of $\Gamma_S$ embedded in $G_T$ for a finite set of places $T$ disjoint from $S$, but we will not consider this extension here to avoid cumbersome notations.

\subsection{Finiteness of the integrability exponent}

Our first task is to show that the integrability exponent defined in \eqref{eq:p_v} is finite:

\begin{thm}\label{th:integrability}
	The integrability exponent $\mathfrak{p}_S({\sf G})$ is finite.
\end{thm}

\begin{proof}
	When ${\sf G}$ is simply connected, this was already established in \cite{GMO}, and we will reduce the proof to this case. We first consider the case when $S$ is finite,
	and then deal with the general case.
	
	We recall that when ${\sf G}$ is simply connected,
	there is no non-trivial automorphic characters (see, for instance, \cite[Lem.~4.1]{GGN1}),
	and it follows from \cite[Th.~3.20 and Th.~3.7]{GMO}
	that for all functions $\phi,\psi\in L_0^2({\sf G}(\mathbb{A}_K)/{\sf G}(K))$ that are $U_\infty$-finite for a maximal compact
	subgroup $U_\infty$ of $G_\infty$ and $\cW$-invariant for a compact open
	subgroup $\cW$ of $G_f$,
	
		 the functions
	\begin{equation}\label{eq:matrix}
	g\mapsto \left< \pi^{\rm aut}_S(g)\phi,\psi\right>=\int_{{\sf G}(\mathbb{A}_K)/{\sf G}(K)} \phi (g^{-1}x)\overline{\psi(x)}\, d\mu(x)
	\end{equation}
	are in $L^p(G_S)$ for a uniform $p>1$.

	Now suppose that ${\sf G}$ is not necessarily simply connected,
	and $S$ is finite.
	We fix a maximal compact subgroup $U_\infty$ of $G_\infty$
	and a compact open subgroup $\cW$ of $G_f$.
	We shall show that for any 
	compactly supported $\phi,\psi\in \mathcal{H}^{00}_{\sf G}$ 
	that are $U_\infty$-finite and $\cW$-invariant,
	the matrix coefficients \eqref{eq:matrix}
	are in $L^p(G_S)$. Since the span of such functions is dense in 
	$\mathcal{H}^{00}_{\sf G}$ as we vary $\cW$, this will imply that 
	the integrability exponent is finite.
	
	We consider the simply connected cover ${\sf p}:\tilde {\sf G}\to {\sf G}$
	that induces the map
	$$
	{\sf p}: \tilde Y:= \tilde {\sf G}(\mathbb{A}_K)/\tilde {\sf G}(K)
	\longrightarrow Y:= {\sf G}(\mathbb{A}_K)/{\sf G}(K).
	$$
	We denote by $\tilde \mu$ and $\mu$ the invariant probability measures on the spaces
	$\tilde Y$ and $Y$ respectively. According to \cite[Ch.~8, Prop.~8.8]{PlaRa},
	there is an exact sequence 
	$$
	\tilde {\sf G}(\mathbb{A}_K) \stackrel{\sf p}{\longrightarrow} {\sf G}(\mathbb{A}_K)\longrightarrow \prod_{v\in V_K} H^1(K_v, Z({\sf G})).
	$$
	so that ${\sf p}(\tilde {\sf G}(\mathbb{A}_K))$ is a normal co-abelian subgroup of ${\sf G}(\mathbb{A}_K)$. 
	Let us consider the group $L:=\cW{\sf p}(\tilde {\sf G}(\mathbb{A}_K)){\sf G}(K)$. Clearly, it is also a normal co-abelian subgroup
	of ${\sf G}(\mathbb{A}_K)$. Furthermore, $L$ is open in ${\sf G}(\mathbb{A}_K)$.
	We also consider a subset $Y_0:=L/{\sf G}(K)$ of $Y$.
	Then 
	$$
	Y=\bigsqcup_{\gamma\in \Delta} \gamma Y_0,
	$$
	where $\Delta$ is a set of coset representatives for $L$ in ${\sf G}(\mathbb{A}_K)$.
	Since $Y_0$ is open, $\mu(Y_0)>0$. Thus, using that $\mu(Y)<\infty$,
	we conclude that $\Delta$ is finite. (In fact, this can be also deduced from the finiteness of class number of $\sf G$.)
	
	Let $m_{\cW}$ be the Haar probability measure on $\cW$.
	The unique $L$-invariant probability measure $\mu_0$ on $Y_0$ can be given as
	\begin{equation}\label{eq:mmm1}
	\int_{Y_0} f\, d\mu_0=\int_{w\in \cW}\int_{\tilde{y}\in \tilde Y} f(w{\sf p}(\tilde{y}))\, d\tilde \mu(\tilde{y})dm_{\cW}(w) \quad \hbox{ for $f\in C_c(Y_0)$.} 
	\end{equation}

	Indeed, the invariance of this measure is easy to check using that 
	${\sf p}(\tilde {\sf G}(\mathbb{A}_K))$ is co-abelian, ${\sf p}$ is equivariant, and $\tilde{\mu}$ is invariant. 	On the other hand,  it also follows from invariance that
	\begin{equation}\label{eq:mmm2}
	\mu_0=|{\sf G}(\mathbb{A}_K):L|\cdot \mu|_{Y_0}.
	\end{equation}
	
	
	Similarly, considering the exact sequence
	$$
	\tilde {G}_S \stackrel{\sf p}{\longrightarrow} {G}_S\longrightarrow \prod_{v\in S} H^1(K_v, Z({\sf G})),
	$$
	we deduce that $H:={\sf p}(\tilde G_S)$ is an open, normal, co-Abelian subgroup of $G_S$. Moreover, since $S$ is finite, $H$ has finite index in $G_S$.
	We have the decomposition 
	$$
	G_S =\bigsqcup_{\delta\in \Omega} \delta H,
	$$
	where $\Omega$ is a finite subset of coset representatives for $H$ in $G_S$. Let $\tilde m_S$ be a Haar
	measure on $\tilde G_S$. 
	Since $H$ is open in $G_S$, a Haar measure on $H$ is given by the restriction of the Haar measure $m_S$ on $G_S$. It follows from the uniqueness of Haar measure that
	for some $c>0$, 
	\begin{equation}\label{eq:mmm3}
	\int_H f\, dm_S=c\int_{\tilde{g}\in \tilde G_S} f({\sf p}(\tilde{g}))\, d\tilde m_S(\tilde{g}) \quad\hbox{for $f\in C_c(H).$} 
	\end{equation}
Indeed, since $S$ is finite, 
	the kernel of the map ${\sf p}:\tilde G_S\to G_S$ is finite,
	so that this map is proper. 
	Thus, the right-hand-side of (\ref{eq:mmm3}) defines a Haar measure on $H$.

	For $\phi,\psi \in \mathcal{H}_{\sf G}^{00}$ satisfying the invariance properties prescribed above, and $g\in G_S$,
	\begin{align*}
	\left< \pi^{\rm aut}_S(g)\phi,\psi\right>=\int_{Y} \phi (g^{-1}y)\overline{\psi(y)}\, d\mu(y)=\sum_{\gamma\in\Delta}
	\int_{Y_0} \phi (g^{-1}\gamma y)\overline{\psi(\gamma y)}\, d\mu(y),
	\end{align*}
	and by H\"older's inequality,
	\begin{align}\label{eq:long}
	&\left(\int_{G_S} \left|\left< \pi^{\rm aut}_S(g)\phi,\psi\right>\right|^p\, dm_S(g)
	\right)^{1/p} \nonumber\\
	\le & \sum_{\gamma\in\Delta}
	\left(\int_{G_S}\left| \int_{Y_0} \phi (g^{-1} \gamma y)\overline{\psi(\gamma y)}\, d\mu(y)\right|^p dm_S(g)\right)^{1/p} \nonumber\\
	\le & \sum_{\gamma\in\Delta}
	\left(\sum_{\delta\in\Omega} \int_H\left| \int_{Y_0} \phi ( h^{-1} \delta^{-1} \gamma y)\overline{\psi(\gamma y)}\, d\mu(y)\right|^p dm_S(h)\right)^{1/p} \nonumber \\
	= & \sum_{\gamma\in\Delta}
	\left(\sum_{\delta\in\Omega} \int_H\left| \int_{Y_0} \phi (\delta^{-1} \gamma h^{-1} y)\overline{\psi(\gamma y)}\, d\mu(y)\right|^p dm_S(h)\right)^{1/p}.
	\end{align}
	Here in the last step we have used the fact that $L$ is co-Abelian. Indeed, denoting $ \delta^{-1}\gamma=z$, this property imples $h^{-1}z=zh^{-1} l'$ with $l'\in L$. Now $Y_0=L/{\sf G}(K)$ and writing a coset $y\in Y_0$ as $y=l {\sf G}(K)$ with $l\in L$, we have $h^{-1}zl {\sf G}(K)=zh^{-1}l'l{\sf G}(K)$. Therefore the integrals over $Y_0=L/{\sf G}(K)$ in the last step above gives the same function of $h\in H$ in both cases. 
	
	We claim that for every $r_1,r_2\in {\sf G}(\mathbb{A}_K)$, the functions
	$$
	c_{r_1,r_2}(h):= \int_{Y_0} \phi (r_1 h y)\overline{\psi(r_2 y)}\, d\mu(y)
	$$
	are in $L^p(H)$. In view of \eqref{eq:mmm1} and \eqref{eq:mmm2},
	\begin{align*}
	c_{r_1,r_2}(h)=|{\sf G}(\mathbb{A}_K):L|^{-1}
	\int_{\tilde{y}\in \tilde Y} \phi (r_1 h {\sf p}(\tilde{y}))
	\overline{\psi(r_2 {\sf p}(\tilde{y}))}\, d\tilde\mu(\tilde{y}),
	\end{align*}
		For a function $f$ on 
	$Y$ and $g\in {\sf G}(\mathbb{A}_K)$, we define a function on $\tilde Y$ by
	$$
	\tilde f_g(\tilde{y}):=f(g {\sf p}(\tilde{y}))\quad \hbox{for $\tilde{y}\in \tilde Y.$} 
	$$
	Then when $f$ is $\cW$-invariant, using \eqref{eq:mmm1}, \eqref{eq:mmm2},
	and that $L$ is co-Abelian, we deduce that 
	\begin{align*}
	\int_{\tilde Y} \tilde f_g\, d\tilde \mu=&\int_{\cW}\int_{\tilde Y} f(g w{\sf p}(\tilde{y}))\, d\tilde \mu(\tilde{y})dm_{\cW}(w) 
	=\int_{\cW}\int_{\tilde Y} f(wg {\sf p}(\tilde{y}))\, d\tilde \mu(\tilde{y})dm_{\cW}(w)\\ 
	=&\int_{Y_0} f(gz)\, d\mu_0(z)
	=|{\sf G}(\mathbb{A}_K):L|\int_{g Y_0} f\, d \mu.
	\end{align*}
	Let $\mathcal{X}({\sf G},L)$ denote the set of continuous unitary characters $\chi$ of ${\sf G}(\mathbb{A}_K)$ such that $\chi(L)=1$.
	It follows from the properties of characters of finite abelian groups that
	$$
	\chi_{Y_0} =|{\sf G}(\mathbb{A}_K):L|^{-1} \sum_{\chi\in\mathcal{X}({\sf G},L)}\chi,
	$$ 
	so that if $f\in \mathcal{H}_{\sf G}^{00}$, using the previous two identities we deduce 
	\begin{align*}
	\int_{\tilde Y} \tilde f_g\, d\tilde \mu&= \int_Y
	(f\circ g)\left({\sum}_{\chi\in\mathcal{X}({\sf G},L)}\chi \right)\, d \mu=
	\int_{Y}  f\left({\sum}_{\chi\in\mathcal{X}({\sf G},L)}\chi\circ g^{-1} \right)\, d \mu \\
	&=
	\sum_{\chi\in\mathcal{X}({\sf G},L)} \chi(g)^{-1}\int_{Y} f\chi\, d \mu=0.
	\end{align*}
	This shows that if $f\in \mathcal{H}_{\sf G}^{00}$, then $\tilde f_g\in L^2_0(\tilde Y)$.
	In particular, $\tilde \phi_{r_1}, \tilde\psi_{r_2}\in L^2_0(\tilde Y)$.
	Using H\"older's inequality, the fact that $L$ is co-Abelian, and \eqref{eq:mmm3}, we conclude that
	\begin{align*}
	\int_H |c_{r_1,r_2}|^p \, dm_S= &
	\int_H \left|\int_{\cW}\int_{\tilde Y} \phi (r_1 h w{\sf p}(\tilde{y}))\overline{\psi(r_2 w{\sf p}(\tilde{y}))}\, d\tilde\mu(\tilde{y})dm_{\cW}(w)\right|^p \,dm_S(h)
	\\
	\le &
	\int_H \int_{\cW} \left|\int_{\tilde Y} \phi (r_1 h w{\sf p}(\tilde{y}))\overline{\psi(r_2 w{\sf p}(\tilde{y}))}\, d\tilde\mu(\tilde{y})\right|^p \,dm_{\cW}(w)dm_S(h)
	\\
	\ll &
	\int_{\cW}\int_H \left|\int_{\tilde Y} \phi (r_1 wh {\sf p}(\tilde{y}))\overline{\psi(r_2 w{\sf p}(\tilde{y}))}\, d\tilde\mu(\tilde{y})\right|^p \,dm_S(h)dm_{\cW}(w)
	\\
	\ll &
	\int_{\cW}\int_{\tilde G_S} \left|\int_{\tilde Y} \phi (r_1 w{\sf p}(\tilde{g}) {\sf p}(\tilde{y}))\overline{\psi(r_2 w{\sf p}(\tilde{y}))}\, d\tilde\mu(\tilde{y})\right|^p d\tilde m_S(\tilde{g})dm_{\cW}(w)\\
	= & \int_{\cW}\int_{\tilde G_S} \left|\left<\tilde{\pi}^{\rm aut}_S(\tilde{g}^{-1}) \tilde \phi_{r_1w}, \tilde\phi_{r_2w}\right> \right|^p d\tilde m_S(\tilde{g})dm_{\cW}(w).
	\end{align*}
	It is easy to see that $\tilde \phi_{r_1w}$ and $\tilde\phi_{r_2w}$
	are $\widetilde U_\infty$-finite for suitable maximal compact subgroup $\widetilde U_\infty$ of $\widetilde{G}_\infty$
	and $\widetilde \cW$-finite for suitable compact open subgroup $\widetilde \cW$ of $\tilde G$. 
	We observe that the integrand in the above formula
	is locally constant in $w\in \cW$, and
	it follows from the simply connected case that 
	the functions $g\mapsto \left<\tilde{\pi}^{\rm aut}_S(\tilde{g}^{-1}) \tilde \phi_{r_1w}, \tilde\phi_{r_2w}\right>$ are in $L^p(\tilde G_S)$.
	Hence, we conclude that
	$$
	\int_{\cW}\int_H |c_{r_1w,r_2w}|^p \, dm_Sdm_{\cW}(w)<\infty,
	$$
	and by \eqref{eq:long} also
	$$
	\int_{G_S} \left|\left< \pi^{\rm aut}_S(g)\phi,\psi\right>\right|^p\, dm_S(g)<\infty.
	$$
	This proves that the representation $\pi^{\rm aut}_S|_{\mathcal{H}_{\sf G}^{00}}$ is $L^p$-integrable when $S$ is finite, and $p$ is uniform and independent of $S$.
	
	\vspace{0.2cm}
	
	To deal with the general case with $S$ arbitrary and possibly infinite, we use that the representation $\pi^{\rm aut}_S$
	has the direct integral decomposition
	\begin{equation}\label{eq:dd1}
	\pi^{\rm aut}_S|_{\mathcal{H}_{\sf G}^{00}} =\int_{\widehat G_S} \rho^{\oplus n(\rho)} \, d\Pi(\rho)
	\end{equation}
	with respect to a measure $\Pi$ on the unitary dual $\widehat G_S$, with $n(\rho)\in \NN\cup\set{\infty}$ denoting the multiplicity.
	The irreducible representations $\rho$ in this decomposition are restricted tensor products of the form
	$\rho=\otimes_{v\in S} \rho_v,$
	where $\rho_v$ are irreducible representations of $G_v$, which are $U_v$-spherical for almost all $v\in S$ (see \cite{f}). For each $v\in S$,
	we also have the decomposition 
	\begin{equation}\label{eq:dd2}
	\pi^{\rm aut}_v|_{\mathcal{H}_{\sf G}^{00}} =\int_{\widehat G_v} \rho_v^{\oplus n_v(\rho_v)} \, d\Pi_v(\rho_v).
	\end{equation}
	We can conclude from the discussion in the previous paragraph that the representation $\pi^{\rm aut}_v|_{\mathcal{H}_{\sf G}^{00}}$ is $L^p$-integrable. This condition implies that $\Pi_v$-almost all $\rho_v$ is $L^{p'}$-integrable for some finite $p' \ge p$. Indeed, $L^p$-integrability implies that a suitable tensor power $\left(\pi^{\rm aut}_v|_{\mathcal{H}_{\sf G}^{00}}\right)^{\otimes N}$ is weakly contained in the regular representation of $G_v$. Since $\rho_v$ is weakly contained in $\pi^{\rm aut}_v|_{\mathcal{H}_{\sf G}^{00}}$, it follows that for $\Pi_v$-almost all $\rho_v$,  the tensor power $\rho_v^{\otimes N}$ is weakly 
	contained in the regular representation of $G_v$.   
	 It follows that the $N$-th power of the matrix coefficients of $U_v$-finite vectors $\rho_v$ satisfy the pointwise bound given by \cite[Thm. 2]{CHH}. Therefore $\rho_v$ is $p^\prime$-integrable for $\Pi_v$-almost all $\rho_v$, and the matrix coefficients of $U_v$-finite vectors  $\phi_v,\psi_v$ of $\rho_v$ satisfy, for suitable $k\in \N$, 		 \begin{equation}\label{eq:bb1}
	\big|\left\langle \rho_v(g_v)\phi_v,\psi_v\right\rangle\big| \le d_v(\phi_v)^{1/2} d_v(\psi_v)^{1/2}\|\phi_v\| \|\psi_v\|\, \Xi_v(g_v)^{1/k}\quad\hbox{for $g_v\in G_v$},
	\end{equation}
	where $d_v(\phi_v):=\dim \left\langle \rho_v(U_v) \phi_v \right\rangle$, and  $\Xi_v$ denotes the Harish-Chandra function on $G_v$. 
	Since the measure $\Pi_v$ in \eqref{eq:dd2} is the image of the measure $\Pi$ 
	from \eqref{eq:dd2} under the restriction map,
	$\Pi$-almost every representation $\rho$ appearing in the decomposition \eqref{eq:dd1} is of the form $\rho=\otimes_{v\in S} \rho_v$
	where the representations $\rho_v$ satisfy the bound \eqref{eq:bb1}.
	It follows from the description of the space of $U_S$-finite vectors for tensor products (see \cite{f}), that for all $U_S$-finite vectors $\phi$ and $\psi$, there exists
	$c(\phi,\psi)>0$ such that 
	\begin{equation}\label{eq:bb2}
	\big|\left\langle \rho(g)\phi,\psi\right\rangle\big| \le c(\phi,\psi)\,\Xi_S(g)^{1/k}
	\quad \hbox{for $g\in G_S$,}
	\end{equation}
	where $\Xi_S(g):=\prod_{v\in S} \Xi_v(g_v)$ is the Harish-Chandra function on $G_S$.
	We recall that the Harish-Chandra function $\Xi_S$ is $L^{4+\eta}$-integrable for all $\eta>0$ 	(see \cite[Prop.~6.3]{GN12}).
	We note that the argument in \cite{GN12} utilizes only the Cartan and Iwasawa decompositions of the group, and does not require that the group be   simply connected. Hence, the estimate \eqref{eq:bb2} implies that 
	$\Pi$-almost every representation $\rho$ appearing in the decomposition \eqref{eq:dd1} is $L^q$-integrable with a uniform $q$, and with the $L^q$-norm uniformly bounded. We conclude that the representation $\pi^{\rm aut}_S|_{\mathcal{H}_{\sf G}^{00}}$ is $L^q$-integrable.
\end{proof}

\subsection{Mean Ergodic Theorem for general groups} \label{sec:mean2}
We note that Theorem \ref{th:mean_sc2} fails 
if $\sf G$ is not simply connected and the corresponding action of $G_S$
is not even ergodic. Nonetheless, it turns out that an analogue of this estimate holds
if we consider actions on a smaller space.

Let $\mathcal{X}({\sf G},\cI_f)$ denote the set of continuous unitary characters $\chi$ of ${\sf G}(\mathbb{A}_K)$ such that $\chi(\cI_f)=\chi(\Gamma)=1$. 
This set is known to be  a finite abelian
group, and its kernel
$$
G^{\ker}:=\{g\in {\sf G}(\mathbb{A}_K):\, \chi(g)=1\quad\hbox{for all $\xi\in \mathcal{X}({\sf G},\cI_f)$}\}.
$$
is a finite index subgroup
of ${\sf G}(\mathbb{A}_K)$ (see \cite[Lem.~4.4]{GGN1}).
Let $G_\infty^{0}$ denote the connected component of identity in $G_\infty$. 
Since $G_\infty^0$ is a connected semisimple Lie group,
it is clear that $G_\infty^0\subset G_\infty\cap G^{\ker}$.
We also set 
$$
G_S^{\ker}:=G_S\cap G^{\ker}.
$$
We note that $G_S^{\ker}$ is a finite index closed (and open) subgroup of $G_S$.
Let
$$
\Gamma_S^{\ker}:=\Gamma_S\cap (G_\infty^0\times G_S).
$$

\begin{lem}\label{l:con}
	$\Gamma_S^{\ker}\subset G_\infty^0\times G^{\ker}_S$.
\end{lem}

\begin{proof}
	We consider $\Gamma^{\ker}_S$ as 
	a subgroup of  $G_\infty^0\times G_S\times \cI^S\subset {\sf G}(\mathbb{A}_K)$ embedded diagonally. Then for every $\chi\in \mathcal{X}(G,\cI_f)$ and $\gamma\in \Gamma_S^{\ker}$, we have
	$\chi(\gamma,\gamma,\gamma)=1$
	and $\chi(G^0_\infty)=\chi(\cI^S)=1$,
	so that it also follows that $\chi(e,\gamma,e)=1$.
	Hence, $\Gamma_S^{\ker}$ is in fact contained in 
	$G_\infty^0\times G_S^{\ker}\times \cI^S$.
\end{proof}

We consider the space
$$
Y_S^{\ker}:=(G_\infty^0\times G^{\ker}_S)/\Gamma_S^{\ker}
$$
equipped with the unique invariant probability measure $\mu^{\rm ker}_S$,
and the corresponding unitary representations $\pi_S$ of $G^{\ker}_S$
on $L^2(Y_S^{\ker})$. 

\begin{thm} \label{th:mean_general}
	Let $\beta$ be a Haar-uniform 
	probability measure supported on $\cI_S$-bi-invariant bounded subset $B$ of $G^{\ker}_S$. Then there exists $\mathfrak{k}_S({\sf G})>0$ such that
	for every given $\eta > 0$, for  $\phi\in L^2(Y^{\ker}_{S})$, 
	$$
	\left\| \pi_{S}(\beta)\phi-\int_{Y^{\ker}_{S}}\phi\,d\mu^{\rm ker}_{S}\right\|_{L^2(Y^{\ker}_{S})}\ll_{S,\eta}
	m_{S}(B)^{-\mathfrak{k}_{ S}({\sf G})+\eta}\, \norm{\phi}_{L^2(Y^{\ker}_{S})}.
	$$
\end{thm}

\begin{proof}
 Part I.  We first consider the case of $S=V_K^f$, the full set of finite places.
 
  To simplify notations, we set
	$$
	G:={\sf G}(\mathbb{A}_K),\quad \Gamma={\sf G}(K),\quad
	X:=G/\Gamma,\quad X^{\ker}:=G^{\ker}/\Gamma.
	$$
	We equip $X$ and $X^{\ker}$ with the Haar probability measures $\mu$ and $\mu^{\rm ker}$
	respectively. Then by invariance,
	\begin{equation}\label{eq:mu}
	\mu^{\rm ker}=\mu(X^{\ker})^{-1}\cdot \mu|_{X^{\ker}}=|G:G^{\ker}|\cdot \mu|_{X^{\ker}}.
	\end{equation}
	
	Let $B^{\ker}$ denote a bounded $\cI_f$-bi-invariant subset of $G_f\cap G^{\ker}$,
	and $\beta^{\ker}$  is the Haar-uniform probability measure.
	We will first prove an effective mean ergodic theorem for the operators
	$\pi_{V_K^f}^{\rm aut}(\beta^{\ker})$ acting on $L^2(X)$
	for test functions $\psi$ with $\hbox{supp}(\psi)\subset X^{\ker}$.
	
	Our argument proceeds similarly to \cite[Th.~4.5]{GGN1}.
	We observe that the space $L^2(X)$ has the decomposition
	$$
	L^2(X)=\mathcal{H}_{\sf G}^{\rm char}\oplus \mathcal{H}_{\sf G}^{00},
	$$
	where $\mathcal{H}_{\sf G}^{\rm char}$ is the closure of the span of automorphic characters, and $\mathcal{H}_{\sf G}^{00}$ is its orthogonal
	complement. We also write 
	$$
	\mathcal{H}_{\sf G}^{\rm char}=\mathcal{H}\oplus \mathcal{H}',
	$$
	where $\mathcal{H}$ is the finite-dimensional 
	space spanned by $\mathcal{X}(G,\cI_f)$, the automorphic characters and  
	$\mathcal{H}'$ its orthogonal complement.
	Since the measure $\beta^{\ker}$ is $\cI_f$-bi-invariant,
	for all automorphic characters $\chi$ and $u\in \cI_f$,
	$$
	\pi^{\rm aut}_{V_K^f}(\beta^{\ker})\chi=\pi^{\rm aut}_{V_K^f}(\beta^{\ker}\ast \delta_u)\chi=\chi(u^{-1})\pi^{\rm aut}_{V_K^f}(\beta^{\ker})\chi.
	$$
	In particular, it follows that $\pi^{\rm aut}_{V_K^f}(\beta^{\ker})\chi=0$
	when $\chi$ is not $\cI_f$-invariant, and 
	$$
	\pi^{\rm aut}_{V_K^f}(\beta^{\ker})|_{\mathcal{H}'}=0.
	$$
	It follows from Theorem \ref{th:norm2} with $S=V_K^f$ that
	$$
	\left\|\pi_{V_K^f}^{\rm aut}(\beta^{\ker})|_{\mathcal{H}_{\sf G}^{00}}\right\|\ll_{S,\eta} m_S(B)^{-\mathfrak{k}_S({\sf G})+\eta}
	$$
	for some explicit $\mathfrak{k}_S({\sf G})>0$. Hence, we conclude that
	for every $\psi\in L^2(X)$,
	$$
	\left\|\pi^{\rm aut}_{V_K^f}(\beta^{\ker})\psi - 
	\pi^{\rm aut}_{V_K^f}(\beta^{\ker})P_{\mathcal H}\psi\right\|_{L^2(X)}
	\ll_{S} m_S(B)^{-\mathfrak{k}_S({\sf G})}\, \|\psi\|_{L^2(X)},
	$$
	where $P_{\mathcal H}$ denotes the orthogonal projection on the space ${\mathcal{H}}$. 
	
	Using the fact that  
	$\mathcal{X}({\sf G},\cI_f)$
	forms an orthonormal basis of $\mathcal{H}$, under our additional assumption that  $\supp(\psi)$  is contained in the intersction of the  kernels of
	$\chi\in \mathcal{X}({\sf G},\cI_f)$, we obtain that 
	\begin{align*}
	P_{\mathcal{H}}\psi=\sum_{\chi\in \mathcal{X}({\sf  G},\cI_f)} \inn{\psi,\chi}_{L^2(X)}\chi
	= \left( \int_X \psi\, d\mu\right)
	\sum_{\chi\in \mathcal{X}({\sf  G},\cI_f)}\chi
	= \left(\int_{X} \psi\, d\mu\right)
	\xi,
	\end{align*}
	where
	$$
	\xi(g):={\sum}_{\chi\in \mathcal{X}({\sf G},\cI_f)} \chi(g)\quad\hbox{for $ g\in G$.}
	$$
	Since $\supp(\beta^{\ker})\subset G^{\ker}$, it follows that
	$$
	\pi^{\rm aut}_{V_K^f}(\beta^{\ker})\chi=\chi\quad\hbox{for $\chi\in \mathcal{X}({\sf G},\cI_f)$}.
	$$
	Hence, we conclude that for all $\psi\in L^2(X^{\ker})$,  
	$$
	\left\|\pi^{\rm aut}_{V_K^f}(\beta^{\ker})\psi - 
	\left(\int_{X} \psi\, d\mu\right)
	\xi\right\|_{L^2(X)}
	\ll_{S} m_S(B)^{-\mathfrak{k}_S({\sf G})}\, \|\psi\|_{L^2(X)}.
	$$
	Since 
	$$
	\xi|_{X^{\ker}}=|\mathcal{X}({\sf  G},\cI_f)|=|G:G^{\ker}|,
	$$
	we also deduce,
	using \eqref{eq:mu}, that 
	for all $\psi\in L^2(X^{\ker})$,
	\begin{equation}\label{eq:mean0}
	\left\|\pi^{\rm aut}_{V_K^f}(\beta^{\ker})\psi - 
	\int_{X^{\ker}} \psi\, d\mu^{\rm ker}\right\|_{L^2(X^{\ker})}
	\ll_{S} m_S(B)^{-\mathfrak{k}_S({\sf G})}\, \|\psi\|_{L^2(X^{\ker})}.
	\end{equation}

Part II.  We now consider a general set of places $S\subset V_f$ with $\sf G$ isotropic over $S$. 	Let $Z_S$ denote the orbit of $G_\infty^0\times G^{\ker}_S\times \cI^S\subset G^{\ker}$
	acting on the identity coset in the space $X^{\ker}$.
	It is open and closed subset of $X^{\ker}$.
	We equip $Z_S$ with the probability measure $\nu_S:=\mu^{\rm ker}(Z_S)^{-1}\mu^{\rm ker}.$ 
	We have an isomorphism 
	$$
	Y_S^{\ker}\simeq \cI^S\backslash Z_S
	$$
	of $(G_\infty^0\times G^{\ker}_S)$-spaces. Therefore, the unitary representation
	$\pi_S|_{G^{\ker}_S}$ on $L^2(Y_S^{\ker})$ is equivalent to the unitary representation
	$G^{\ker}_S$ on  $L^2(\cI^S\backslash Z_S)$, which is also equivalent 
	to the unitary representation of $G^{\ker}_S$ on 
	the space $L^2(Z_S)^{\cI^S}$ consisting of 
	$\cI^S$-invariant functions in $L^2(Z_S)$.
	More explicitly, given a function $\phi$ on the space $Y^{\ker}_S$, we get
	a $\cI^S$-invariant function
	$$
	\phi^{\ker}\big((g_\infty,g_S,u)\Gamma\big):=\phi\big((g_\infty,g_S)\Gamma_S^{\ker}\big),\quad \hbox{ for $(g_\infty,g_S,u)\in G_S^{\ker}\times G^{\ker}_S\times \cI^S,$}
	$$
	on the space $Z_S$. This defines a $(G_\infty^0\times G^{\ker}_S)$-equivariant
	isometry between $L^2(Y^{\ker}_S)$ and $L^2(Z_S)^{\cI^S}$.
	
	 Recall that we denote by $\beta$ the Haar-uniform 
	probability measure supported on $\cI_S$-bi-invariant bounded subset $B$ of $G^{\ker}_S$.
	Denote by  $\beta^{\ker}$ the Haar-uniform probability measure
	supported on $B\times \cI^S\subset G_f$. Then since $\phi^{\ker}$ is $\cI^S$-invariant,
	we obtain that 
	$$
	\pi^{\rm aut}_{V_K^f}(\beta^{\ker})\phi^{\ker}=\pi_S(\beta)\phi.
	$$
	Finally, we apply the estimate \eqref{eq:mean0} to the case when $\psi=\phi^{\ker}$
	for $\phi\in L^2(Y^{\ker}_S)$. Since 
	$$
	\int_{X^{\ker}} \phi^{\ker}\, d\mu^{\rm ker}=\mu^{\rm ker}(Z_S)\int_{Z_S}\phi^{\ker}\, d\nu_S
	=\mu^{\rm ker}(Z_S)\int_{Y^{\ker}_S}\phi\, d\mu^{\rm ker}_S,
	$$
	using the above identifications, we conclude that (for any $\eta > 0$)
	$$
	\left\|\pi_{S}(\beta)\phi - \mu^{\rm ker}(Z_S)
	\int_{Y^{\ker}_S} \phi\, d\mu^{\rm ker}_S
	\right\|_{L^2(Y^{\ker}_S)}
	\ll_{S,\eta} m_S(B)^{-\mathfrak{k}_S({\sf G})+\eta}\, \|\phi\|_{L^2(Y^{\ker}_S)}.
	$$
	We note that this estimate holds for any bounded $\cI_S$-bi-invariant
	subsets of $G_S^{\ker}$, so that taking $\phi=1$ and $m_S(B)\to\infty$,
	we deduce from the above estimate that $\mu^{\ker}(Z_S)=1$.
	This completes the proof of the Theorem \ref{th:mean_general}.
\end{proof}

\subsection{Discrepancy estimates for general groups}
When $\sf G$ is not simply connected, the set $\Gamma_S$ may be not dense in $G_\infty$. Nonetheless, according to Corollary \ref{c:dense}, its closure is a subgroup of finite index in $G_\infty$. In particular, $\Gamma_S$ is dense in $G^0_\infty$, the connected component of $G_\infty$, and we estimate the discrepancy for $\Gamma_S$-points for subsets of $G^0_\infty$.
As in the previous results, we parametrize  $\Gamma_S$ by the subsets
$$
R_S(h):=\{\gamma\in \Gamma_S:\, \h(\gamma)\le h \}.
$$
Let $m^0_\infty$ be the Haar measure on $G_\infty^0$, which we choose to  normalize so that
$\Gamma_S$ has covolume one in $G_\infty^0\times G_S^{\ker}$ with respect to $m^0_\infty\times m_S$.
We define 
$$
v^{\ker}_S(h):=m_S(B^{\ker}_h(S)),
$$
where $B^{\ker}_h(S):=\{g\in G^{\ker}_S:\, \h(g)\le h \}$.
For $\Omega\subset G^0_\infty$, we introduce the discrepancy function:
$$
\mathcal{D}(R_S(h),\Omega):= \left|\frac{|R_S(h)\cap \Omega|}{v^{\ker}_S(h)}-m^0_\infty(\Omega)\right|.
$$
We emphasize here that the correct normalisation is in fact by $v^{\ker}_S(h)$, and not by $v_S(h)$ as in the simply-connected case. While $v^{\ker}_S(h)$ is comparable with  $v_S(h)=m_S(B_h(S))$
up to multiplicative constants, to get the correct main term it is essential to  take into account the contribution of automorphic characters. This contribution manifests itself through the volume function $v^{\ker}_S(h)$ of $B_h^{\ker}\subset G_S^{\ker}$. Note that looking only at the group of  rational points $\Gamma_S$ which is a lattice in $G_\infty\times G_S$ as before, the fact $B_S^{\ker}(h)$ is the correct choice here is not obvious  in advance. It reflects subtle algebraic information regarding the behavior of the automorphic characters when restricted to $G_S$, since in fact  $\Gamma_S\subset G_\infty^{0}\times G_S^{\ker}$. 

Using the methods presented in \S 4, it  is possible to establish discrepancy results in the present case given congruence constraints determined by suitable compact-open subgroups $W\subset \cI^S$. But since in the non-simply-connected case certain congruence obstructions are bound to arise, the full analysis here is longer and requires additional notation. For brevity, we state the results only for the case of $W=\cI^S$, namely in the absence of congruence constraints.  

As above for $x\in G^0_\infty$ and $E\subset G^0_\infty$, we set $E(x):=Ex$.

\medskip

Using Theorem \ref{th:mean_general}, we establish 
a mean-square discrepancy bound:

\begin{thm}[mean-square discrepancy bound] \label{th:l2_bound_general}
	Let $E$ be any measurable subset of $G^0_\infty$ of positive finite measure 
	satisfying $\mathcal{N}(E)<\infty.$
Then for any $\eta > 0$ 
	$$
	\Big\|\mathcal{D}(R_S(h), E(\cdot))\Big\|_{L^2(Q)}\ll_{S, Q, \eta} \mathcal{N}(E)^{1/2}
	m^0_\infty(E)^{1/2} v^{\ker}_S(h)^{-\mathfrak{k}_S({\sf G})+\eta}
	$$ 
for every bounded measurable subset $Q$ of $G^0_\infty$.	
\end{thm}

\begin{proof}
The argument proceeds along the lines of the proof of Theorem \ref{th:l2_bound}.
We note that since $E(x)$ is a subset of $G^0_\infty$, by Lemma \ref{l:con},
$$
\Gamma_S\cap (E(x)\times B_S(h))=\Gamma_S\cap (E(x)\times B^{\ker}_S(h))=\Gamma^{\ker}_S\cap (E(x)\times B^{\ker}_S(h)).
$$
We consider the function 
$$
\phi(g):=\sum_{\delta\in \Gamma^{\ker}_S} \chi_{E^{-1}}(g_\infty\delta^{-1}) \chi_{\cI_S} (g_S\delta^{-1})
$$ 
on $Y_S^{\ker}=(G_\infty^0\times G^{\ker}_S)/\Gamma_S^{\ker}$.
As in the proof of Theorem \ref{th:l2_bound}, one verifies that 
$$
\big|\Gamma^{\ker}_S\cap (E(x)\times B^{\ker}_S(h))\big|= \int_{a\in B^{\ker}_S(h)}\phi\big(a(x,e)\big)\, dm_{S}(b),
$$

and 
\begin{align*}
\mathcal{D}(R_S(h), E(x)) &=
\left| \frac{1}{m_S(B^{\ker}_S(h))}\int_{a\in B^{\ker}_S(h)}\phi\big(a(x,e)\big)\, dm_{S}(b)-\int_{Y^{\ker}_{S}} \phi\, d\mu_{S}  \right|.
\end{align*}
We observe that since the subgroup $G^{\ker}$ is defined as a kernel of a set of $\cI_S$-invariant characters, it is $\cI_S$-bi-invariant. Hence, it follows
that the sets $B^{\ker}_S(h)=G^{\ker}\cap B_S(h)$ are also $\cI_S$-bi-invariant,
and we can apply Theorem \ref{th:mean_general}. The remaining proof proceeds exactly 
as the proof of Theorem \ref{th:l2_bound}.
\end{proof}

Once the mean-square discrepancy bound has been established, one 
can also deduce exactly as Section \ref{sec:dis_simp} generalizations of 
the almost-sure discrepancy bound (Theorem \ref{almost-sure-disc}),
the uniform discrepancy bound for right-stable sets (Theorem \ref{th:wr_general}),
and the uniform discrepancy bound for balls of arbitrarily small radius (Theorem \ref{th:balls_general}).
Here the pointwise estimates depend as before on the dimension $\mathfrak{d}:=\dim_\RR (G_\infty)$. We state these results as follows.

\begin{thm}[almost-sure discrepancy bound]\label{almost-sure-disc_2}
	With notation as in Theorem~\ref{th:l2_bound_general},
	for every $\eta>0$ and almost all $x\in G^0_\infty$, 
	$$
	\mathcal{D}\big(R_S(h), E(x)\big)\ll_{S,E,x,\mathfrak{k},\eta}  \big(\log v^{\ker}_S(h)\big)^{3/2+\eta} v^{\ker}_S(h)^{-\mathfrak{k}}.
	$$	
\end{thm}

\begin{thm}[uniform discrepancy bound]\label{th:wr_general_2}
	For every right-stable subsets $E$ of $G_\infty$ of finite measure such that $\mathcal{N}(E_{\epsilon_0}^+)<\infty$, for every $0 < \mathfrak{k}= \mathfrak{k}_S({\sf G})-\eta < \mathfrak{k}_S({\sf G})$ 
and for $x\in G_\infty$, the following pointwise bound for the discrepancy holds :
	$$
	\mathcal{D}(R_S(h), E(x)) \ll_{S,E,x,\eta} v^{\ker}_S(h)^{-2\mathfrak{k}/(\mathfrak{d}+2)}
	$$
	provided that $h\ge h_0(S)$.
Moreover, this estimate is uniform for $x$ ranging in compact subsets of $G_\infty$.
\end{thm}

\begin{thm}[discrepancy bound for balls]\label{th:balls_general_2}
	Let $x\in G_\infty$ and $\ell\in (0,\ell_0)$, for suitable $\ell_0$ (independent of $x$.  For every $0 < \mathfrak{k}= \mathfrak{k}_S({\sf G})-\eta < \mathfrak{k}_S({\sf G})$ : 
	$$
	\big|R_S(h)\cap B(x,\ell)\big|= m_\infty(B(e,\ell))v^{\ker}_S(h)+O_{S,x, \eta}\Big (m_\infty(B(e,\ell))^{\mathfrak{d}/(\mathfrak{d}+2)}v_S(h)^{1-2\mathfrak{k}/(\mathfrak{d}+2)}\Big)
	$$

	provided that $m_\infty(B(e,\ell))^2\gg_\eta v_S(h)^{-2\mathfrak{k}}$. 
Explicity, if the volume growth satisfies $v_S(h)\gg_{S,\eta^\prime} h^{\mathfrak{a}}$, with $0< \mathfrak{a}=\mathfrak{a}_S({\sf G})-\eta^\prime$, then the estimate holds provided the height $h$ satisfies
$h\gg_{S,\eta,\eta^\prime}\ell^{-\mathfrak{d}/\mathfrak{k}\mathfrak{a}} \,.$ 
Moreover, this estimate is uniform for $x$ ranging in compact subsets of $G_\infty$.

	\end{thm}

Since the proofs of these results proceed as in Section \ref{sec:dis_simp},
we omit the details.

\end{document}